\documentclass[leqno,11pt]{amsart}
\usepackage{amsmath,amscd,amsthm,amsxtra}
\usepackage{epsfig,graphics,color,colortbl}
\usepackage{amssymb,latexsym}
\usepackage{mathrsfs}
\usepackage{eucal,upgreek}
\usepackage[poly,all]{xy}
\usepackage{hyperref}

\newtheorem{thm}{\bf Theorem}[section]

\newtheorem{prop}[thm]{\bf Proposition}
\newtheorem{cor}[thm]{\bf Corollary}
\newtheorem{lem}[thm]{\bf Lemma}

\newtheorem{ex}[thm]{\bf Example}

\newenvironment{red}
{\relax\color{red}}
{\hspace*{.5ex}\relax}

\newcommand{\ber}{\begin{red}}
\newcommand{\er}{\end{red}}

\newcommand{\A}{\mathcal{A}}

\newcommand{\B}{\mathcal{B}}

\newcommand{\cP}{\mathscr{P}}

\newcommand{\pf}{\noindent{\bfseries Proof. }}

\newcommand{\F}{\mathcal{F}}

\newcommand{\gl}{\mathfrak{gl}}
\newcommand{\Z}{\mathbb{Z}}

\newcommand{\te}{\widetilde{e}}
\newcommand{\tf}{\widetilde{f}}
\newcommand{\tte}{\widetilde{\mathsf e}}
\newcommand{\ttf}{\widetilde{\mathsf f}}
\newcommand{\g}{\mathfrak{g}}

\newcommand{\mc}{\mathcal}
\newcommand{\mf}{\mathfrak}

\newcommand{\hyphen}{\raisebox{.2ex}{\text{-}}}

\newcommand{\sh}{\mathrm{sh}}
\newcommand{\SST}{SST}
\newcommand{\wcol}{{w}_{\mathrm{col}}}
\newcommand{\wrow}{{w}_{\mathrm{row}}}
\newcommand{\rhocol}{\varrho_{\mathrm{col}}}
\newcommand{\rhorow}{\varrho_{\mathrm{row}}}

\newcommand{\cM}{\mathbf{M}}

\numberwithin{equation}{section}

\begin{document}
\title[Duality on Fock spaces and combinatorial energy functions]
{Duality on Fock spaces and combinatorial energy functions}
\thanks{J.-H. Kwon was supported by Basic Science Research Program through the National Research Foundation of Korea (NRF) funded by the Ministry of  Education, Science and Technology (No. 2011-0006735).}
\author{JAE-HOON KWON}
\address{ Department of Mathematics \\ Sungkyunkwan University \\ Suwon,  Republic of Korea}
\email{jaehoonkw@skku.edu}

\author{EUIYONG PARK}
\thanks{E. Park was supported by the National Research Foundation of Korea(NRF) Grant funded by the Korean Government(MSIP)(NRF-2014R1A1A1002178).}
\address{ Department of Mathematics \\ University of Seoul  \\ Seoul,  Republic of Korea}
\email{epark@uos.ac.kr}

\begin{abstract}
We generalize in a combinatorial way the notion of the energy function of affine type $A$ on a sequence of row or column tableaux to the case of a more general class of modules over a general linear Lie superalgebra $\g$ based on a Howe duality of type $(\g,\gl_n)$ on various Fock spaces.
\end{abstract}

\maketitle

%\setcounter{tocdepth}{1}
%\tableofcontents

\section{Introduction}
The Kostka-Foulkes polynomials are natural $q$-deformation of Kostka numbers. They appear as the entries of a transition matrix between Schur functions and Hall-Littlewood functions, and also coincide with  the Lusztig's $q$-weight multiplicities of type $A$ (cf.~\cite{Mac95}). One of the most important and interesting properties of Kostka-Foulkes polynomials is that they have non-negative integral coefficients. In \cite{LS}  Lascoux and Sch\"{u}tzenberger introduced the notion of charge statistic on semistandard Young tableaux, and proved this positivity of Kostka-Foulkes polynomials.

In \cite{NaYa} Nakayashiki-Yamada showed that the energy function on a finite affine crystal associated to a tensor product of symmetric (or exterior) powers of the natural representation of $\gl_\ell$ is given by the Lascoux and Sch\"{u}tzenberger's charge. This also gives another combinatorial realization of the Kostka-Foulkes polynomials as $q$-deformed decomposition multiplicities.  Indeed, if we understand the energy function on a sequence of row (or column) tableaux as a statistic on the corresponding non-negative integral (or binary) matrix, then the Nakayashiki-Yamada's result implies that the energy of a given matrix
%and hence its corresponding insertion tableau or $P$-tableau
is equal to the charge of its associated recording tableau under the RSK correspondence.
%(Kwon: This part was modified as a response to comment no.1 by the first reviewer)
We also refer the reader to for example, \cite{KSS,SW99,S} and references therein for a generalization of the work \cite{NaYa} to the case of a tensor product of arbitrary Kirillov-Reshetikhin crystals of affine type $A$.
%(Kwon: This part is a response to comment (2) by the 3rd reviewer)

The purpose of this paper is to understand and generalize the energy function of type $A_{\ell-1}^{(1)}$ on a sequence of row (or column) tableaux in another direction from a viewpoint of duality principle due to Howe \cite{H}. Let $(\A,\B)$ be a pair of countable $\Z_2$-graded totally ordered sets. We consider a Lie superalgebra $\g$ of type $A$ and a semisimple $\g$-module $\mathscr{F}$ associated to $(\A,\B)$.
Following \cite{H} (see also \cite{CW12}), one can show that $\g$ forms a dual pair with $\gl_n$ on $\mathscr{F}^{\otimes n}$ for $n\geq 1$, giving a family of irreducible $\g$-modules $L_\g(\lambda)$ appearing in the $(\g,\gl_n)$-decomposition
\begin{equation}\label{Eq:HoweDuality}
\mathscr{F}^{\otimes n}\cong \bigoplus_{\lambda\in H_{\g,n}} L_\g(\lambda)\otimes L_n(\lambda),
\end{equation}
where $H_{\g,n}$ is a subset of $\Z_+^n$, the set of generalized partitions of length $n$, and $L_n(\lambda)$ is the finite-dimensional irreducible $\gl_n$-module corresponding to $\lambda$ (Theorem \ref{Thm:HoweDuality}).

The decomposition \eqref{Eq:HoweDuality} is
the classical $(\gl_\ell,\gl_n)$-duality with $L_\g(\lambda)$ an irreducible polynomial $\gl_\ell$-module, when $\A$ is finite with $\ell$ even elements and $\B=\emptyset$. Moreover, under suitable choices of $(\A,\B)$, $\{\,L_\g(\lambda)\,|\,n\geq 1, \lambda\in H_{\g,n}\,\}$ may also include various interesting families of irreducible modules, which forms a semisimple tensor category, for example, the integrable highest weight modules over $\gl_\infty$, infinite-dimensional unitarizable modules over $\gl_{p+q}$ called a holomorphic discrete series, the irreducible polynomial modules over a general linear Lie superalgebra $\gl_{p|q}$, and so on (cf.~\cite{CL,CLZ,CW03,CW12,Fr,H,KacR2,KV}). A uniform combinatorial character formula for $L_\g(\lambda)$ was given by the first author \cite{K08} in terms of certain pairs of Young tableaux, which we call {\it parabolically semistandard tableaux of shape $\lambda$  (of level $n$)}.

Now, we consider a $\g$-module
\begin{equation}\label{Eq:V(mu)}
V_\g(\mu)=L_\g(\mu_1)\otimes \cdots\otimes L_\g(\mu_n)\subset \mathscr{F}^{\otimes n},
\end{equation}
for $\mu=(\mu_1,\ldots,\mu_n)\in \Z_+^n$ with $\mu_1,\ldots,\mu_n\in H_{\g,1}$, which is semisimple and decomposes into $L_\g(\lambda)$'s for $\lambda\in H_{\g,n}$ with finite multiplicity given by the classical Kostka number $K_{\lambda\mu}$ \cite{K08}. Then one can define a $q$-deformed character of $V_\g(\mu)$ by replacing the multiplicity $K_{\lambda\mu}$ of $L_\g(\lambda)$ in $V_\g(\mu)$ with the corresponding Kostka-Foulkes polynomial $K_{\lambda\mu}(q)$, which can be viewed as a natural $\g$-analogue of modified Hall-Littlewood function.

%From the combinatorial model in \cite{K08}, the character of $V_\g(\mu)$ is equal to the weight generating function of $n$-tuples of parabolically semistandard tableaux of level $1$, say, $\mathcal{F}^n_{\A/\B}$.
As a main result, we introduce a purely combinatorial statistic called a {\it combinatorial energy function} on the $n$-tuples of parabolically semistandard tableaux of level $1$ associated to \eqref{Eq:V(mu)}, which generalizes the usual energy function of type $A_{\ell-1}^{(1)}$ on a sequence of row (or column) tableaux, and also produces the $q$-deformed character of $V_\g(\mu)$ in a bijective way (Theorem \ref{Thm:Q=gradedchar}).
The main ingredient of our proof is an analogue of RSK algorithm for the decomposition \eqref{Eq:HoweDuality} of $\mathscr{F}^{\otimes n}$ as a $(\g,\gl_n)$-module \cite{K08}, which is proved here to be an isomorphism of $\gl_n$-crystals (Theorem \ref{Thm:crystal isom kappa}). Another important one is an intrinsic characterization of the charge statistic on regular $\gl_n$-crystals, which is deduced by combining the result in \cite{NaYa} and a bicrystal structure on the classical RSK correspondence (Theorem \ref{Thm: charge}).

We remark that as in the case of the classical $(\gl_\ell,\gl_n)$-duality it would be very interesting to find a representation theoretic meaning of our combinatorial energy function  in terms of representations of a quantum (super)algebra associated to an affinization of $\g$, especially when $\g=\gl_\infty$ with $L_\g(\lambda)$ the integrable highest weight module, or $\g=\gl_{p|q}$ with $L_\g(\lambda)$ the irreducible  polynomial module  \cite{Her,Zh}.

The paper is organized as follows. In Section 2, we recall the notion of parabolically semistandard tableaux and related results. In Section 3, we review the energy function of affine type $A_{\ell-1}^{(1)}$, and the charge statistic on regular $\gl_n$-crystals together with its new intrinsic characterization. In Section 4, we introduce a combinatorial energy function, and then show that the associated $q$-deformed decomposition multiplicities recover the usual  Kostka-Foulkes polynomial.\vskip 2mm

{\bf Acknowledgement} We would like to thank the referees for careful reading and helpful comments.

\section{Parabolically semistandard tableaux}

\subsection{Semistandard tableaux} \label{Sec: SST}
Let us briefly recall necessary background on semistandard tableaux (cf.\ \cite{Fu,K08}).
Let $\cP$ be the set of partitions, where we often identify a partition with its Young diagram as usual.   Throughout this paper, $\A$ (or $\B$) denotes
 a countable  $\Z_2$-graded set (that is, $\A=\A_0\sqcup\A_1$) with a total order $<$.  By convention,
$\Z_{>0}$, $\Z_{\geq 0}$, and $\Z_{<0}$ denote the set of positive, non-negative, and negative integers with the usual total order and even degree, respectively. For each $\Z_{>0}$, we put $[n] = \{1<2< \cdots < n \}$ and $[-n] = \{-n< \cdots<-2 < -1 \}$ with even degree. We assume that $\A'=\{\,a'\,|\,a\in \A\,\}$ is the set with the total order $a'_1<a_2'$ for $a_1<a_2\in\A$ and the opposite $\Z_2$-grading.

For a skew Young diagram $\lambda/\mu$,
an {\it $\A$-semistandard tableau} $T$ of shape $\lambda/\mu$ is
a filling of $\lambda/\mu$ with entries in $\A$ such that
(1) the entries in rows and columns are weakly increasing from left to right and from top to bottom, respectively,
(2) the entries in $\A_0$ are strictly increasing in each column,
(3) the entries in $\A_1$ are strictly increasing in each row.
Let $\sh(T) $ denote the shape of $T$, and  $\wcol(T)$ (resp.\ $\wrow(T)$) the word with letters in $\A$ obtained by reading the entries column by column (resp.\ row by row) from right to left (resp.\ from bottom to top), and in each column (resp.\ row) from top to bottom (resp.\ from left to right).
Let $\SST_\A(\lambda/\mu)$ be the set of all $\A$-semistandard tableaux of shape $\lambda/\mu$. We set $\cP_\A = \{ \lambda \in \cP \mid \SST_\A(\lambda) \neq \emptyset  \}$. For example, $\cP_n :=\cP_{[n]} =\{ \lambda \in \cP \mid \ell(\lambda) \le n \}$, where $\ell(\lambda)$ is the length of $\lambda$.

Let $P_{\A}=\bigoplus_{ a \in \A}\mathbb{Z}\epsilon_a$ be the free
abelian group with the basis $\{\,\epsilon_a\,|\,a \in \A\,\}$, and let ${\bf x}_{\A}=\{\, x_a \,|\, a \in\A\,\}$ be the set of
commuting formal variables indexed by $\A$.
For $T\in \SST_\A(\lambda/\mu)$,
let ${\rm wt}_{\A}(T)=\sum_{a \in\A}m_a \epsilon_a \in  P_{\A}$ be the weight of $T$, where
$m_a$ is the number of occurrences of $a$ in $T$, and put ${\bf x}_{\A}^{T}=\prod_{a\in\A} x_a^{m_a}$. We define the character of
$\SST_{\A}(\lambda/\mu)$  to be $s_{\lambda/\mu}({\bf x}_{\A})=\sum_{T\in
\SST_{\A}(\lambda/\mu)}{\bf x}_{\A}^{T}$.

For $a\in \A$ and $T \in \SST_\A(\lambda)$, $(T \leftarrow a)$ denotes the tableau obtained by
the  Schensted's column bumping algorithm, and $(a \rightarrow T )$ the tableau obtained by the  row bumping algorithm.
For $T \in \SST_\A(\lambda)$ and $T' \in \SST_\A(\mu)$, we set
$(T \leftarrow T') = (((T \leftarrow c_1) \cdots) \leftarrow c_t )$ and $(T' \rightarrow T) = ( r_t  \rightarrow  ( \cdots (r_1 \rightarrow T ) ) )$,
where $\wcol(T') = c_1\ldots c_t$ and $\wrow(T') = r_1\ldots r_t$.

Let $\mu=(\mu_1, \ldots, \mu_r)$ be a sequence of non-negative integers. For $(T_1, \ldots, T_r) \in \SST_\A(\mu_1) \times \cdots \times \SST_\A(\mu_r)$, let $S_k = (( (T_1 \leftarrow T_2 ) \cdots  ) \leftarrow T_k) $ for $1\leq k\leq r$. We  define $\rhocol(T_1, \ldots, T_r) = (S, S_R)$, where
$S = S_r$ and $S_R$ is the $[r]$-semistandard tableau of shape $\sh(S)$ obtained by filling $\sh(S_{k}) / \sh(S_{k-1})$ with $k$ for $1\leq k\leq r$.
Similarly, we define $\rhorow(T_1, \ldots, T_r) = (S', S_R')$, where $S' = ( T_r \rightarrow  (  \cdots  (T_2 \rightarrow T_1) ) )$. Then we have bijections
\begin{equation}\label{Eq:rho}
 \rhocol, \rhorow : \SST_\A(\mu_1) \times \cdots \times \SST_\A(\mu_r) \longrightarrow \bigsqcup_{\lambda \in \cP_\A} \SST_\A(\lambda) \times \SST_{[r]}(\lambda)_\mu,
 \end{equation}
where $\SST_{[r]}(\lambda)_\mu=\{\,T\in \SST_{[r]}(\lambda)\,|\,{\rm wt}_{[r]}(T)=\sum_{i=1}^r\mu_i\epsilon_i\,\}$.

Let $\A^\pi $ be $\A$ as a $\Z_2$-graded set with the reverse total order of $\A$. For $T\in \SST_{\A}(\lambda/\mu)$, define $T^{\pi}$ to be the tableau obtained after $180^\circ$-rotation of $T$, which is an $\A^\pi$-semistandard tableau.
Let $\A \ast \B$ be the $\Z_2$-graded set $\A \sqcup \B$ with the extended total order defined by $x<y$ for all $x\in \A$ and $y\in \B$.
For $S \in \SST_\A(\mu)$ and $T \in \SST_\B(\lambda/\mu)$, define $S\ast T$ to be the tableau of shape $\lambda$ given by gluing $S$ and $T$ so that $S\ast T \in \SST_{\A\ast \B}(\lambda)$.

\begin{ex} \label{Ex: Ex1}
{\rm
Note that the ordered sets $\Z_{>0}'$ and $\Z_{\ge0}'$ have only elements of odd degree. Letting
$
P = \resizebox{.14\hsize}{!}{${\def\lr#1{\multicolumn{1}{|@{\hspace{.6ex}}c@{\hspace{.6ex}}|}{\raisebox{-.3ex}{$#1$}}}\raisebox{-.0ex}
{$\begin{array}{cccc}
\cline{1-4}
 \lr{1'} & \lr{2'} & \lr{3'} & \lr{4'} \\
\cline{1-4}
\lr{1'}  & \lr{2'} &  \lr{4'} &     \\
\cline{1-3}
\end{array}$}}$}
\quad \in \SST_{\Z_{>0}'}(4,3)
$, we have
$$P^\pi=
\resizebox{.14\hsize}{!}{${\def\lr#1{\multicolumn{1}{|@{\hspace{.6ex}}c@{\hspace{.6ex}}|}{\raisebox{-.3ex}{$#1$}}}\raisebox{-.0ex}
{$\begin{array}{cccc}
\cline{2-4}
   & \lr{4'} & \lr{2'} & \lr{1'} \\
\cline{1-4}
\lr{4'} & \lr{3'}  & \lr{2'} &  \lr{1'}     \\
\cline{1-4}
\end{array}$}}$}\ \quad \in \SST_{\left(\Z_{>0}'\right)^\pi}((4,4)/(1)).
$$
If we consider
\begin{align*}
H &=
\resizebox{.125\hsize}{!}{${\def\lr#1{\multicolumn{1}{|@{\hspace{.6ex}}c@{\hspace{.6ex}}|}{\raisebox{-.3ex}{$#1$}}}\raisebox{-.0ex}
{$\begin{array}{cccc}
\cline{1-4}
\lr{1}  & \lr{1} & \lr{1} & \lr{1} \\
\cline{1-4}
\lr{2} &  \\
\cline{1-1}
\end{array}$}}$}
\  \qquad \qquad \qquad \qquad \in \SST_{[3]}( 4,1) , \\
T &=
\resizebox{.32\hsize}{!}{${\def\lr#1{\multicolumn{1}{|@{\hspace{.6ex}}c@{\hspace{.6ex}}|}{\raisebox{-.3ex}{$#1$}}}\raisebox{-.1ex}
{$\begin{array}{cccccccccc}
& & & & & & & & \\
\cline{5-9}
   &    &    &    &  \lr{2'} & \lr{4'} & \lr{5'} & \lr{6'}  & \lr{7'} \\
\cline{2-9}
   &  \lr{ 0'} & \lr{ 3'} & \lr{ 4'} & \lr{ 5'}   \\
\cline{1-5}
\lr{ 0'} & \lr{ 1'} & & & & & & & \\
\cline{1-2}
& & & & & & & & \\
\end{array}$}}$}
\quad \in \SST_{\Z_{\ge 0}'}((9,5,2)/(4,1)),
\end{align*}
then we have
$$
H*T =
\resizebox{.32\hsize}{!}{${\def\lr#1{\multicolumn{1}{|@{\hspace{.6ex}}c@{\hspace{.6ex}}|}{\raisebox{-.3ex}{$#1$}}}\raisebox{-.1ex}
{$\begin{array}{cccccccccc}
\cline{1-9}
\lr{\,1\,}  & \lr{\,1\,}  & \lr{\,1\,}  & \lr{\,1\,}  &  \lr{2'} & \lr{4'} & \lr{5'} & \lr{6'}  & \lr{7'} \\
\cline{1-9}
 \lr{\,2\,}  &  \lr{0'} & \lr{3'} & \lr{4'} & \lr{5'}   \\
\cline{1-5}
\lr{0'} & \lr{1'} \\
\cline{1-2}
\end{array}$}}$}
\qquad \in \SST_{[3]*\Z_{\ge0}'}(9,5,2).
$$
%\cmt{(Kwon: I removed the part about identifying $(\Z_{>0}')^\pi$ and $\Z_{<0}'$)}
\vskip 0.3em
}
\end{ex}

\subsection{Rational semistandard tableaux}
Let us recall the notion of rational semistandard tableaux \cite{ST87}.
Let $\Z^n_+=\{\,(\lambda_1,\ldots,\lambda_n)\,|\,\lambda_i\in \Z,\, \lambda_1\geq \ldots\geq \lambda_n\,\}$ be the set of generalized partitions of length $n$. We may identify $\lambda\in \Z_+^n$ with a generalized Young diagram. For example, $\lambda=(3,2,0,-2)\in \Z^4_+$ corresponds to
$$\resizebox{.18\hsize}{!}{${\def\lr#1{\multicolumn{1}{|@{\hspace{.6ex}}c@{\hspace{.6ex}}|}{\raisebox{-.3ex}{$#1$}}}\raisebox{-.6ex}
{$\begin{array}{ccc|ccc}
& & & & &  \\
\cline{4-6}
& & & \lr{ } & \lr{ } & \lr{ }  \\
\cline{4-6}
& &   & \lr{ }& \lr{ } & \\
\cline{4-5}
  &   &  &  &  &   \\
\cline{2-3}
  & \lr{ }  &  \lr{  }&  &  & \\
\cline{2-3}
  & \!\!\!{{}_{-2}} & \!\!\!{{}_{-1}} & {{}_{1}} & {{}_{2}} & {{}_{3}}
  \end{array}$}}$}$$\vskip 2mm
\noindent where the non-zero integers indicate the column indices.

For $\lambda \in \Z_+^n$,  a {\it rational semistandard tableau  $T$ of shape $\lambda$} is a filling of $\lambda$ with entries in $[n] \sqcup [-n]$ such that
(1) the subtableau with columns of positive indices is $[n]$-semistandard,
(2)  the subtableau with columns of negative indices is $[-n]$-semistandard,
(3) if $b_1 < \cdots < b_s$ (resp.\ $-b_1' < \cdots < -b_t'$) are the entries in the 1st (resp.\ $-1$st) column with $s+t\leq n$, then $b_i'' \le b_i$ for $1 \le i \le s$, where
$\{ b_1'' < \cdots < b_{n-t}'' \} = [n] \setminus \{ b_1', \ldots , b_t'\}$.
Let us call $n$ the {\it rank of $T$}, and  define the weight of $T$ to be ${\rm wt}_{[n]}(T)=\sum_{i\in [n]}(m^+_i-m^-_i)\epsilon_i$, where $m^\pm_i$ is the number of occurrences of $\pm i$ in $T$.
We also use the same notation $\SST_{[n]}(\lambda)$ to denote the set of rational semistandard tableaux of shape $\lambda$.
%\cmt{(Park: Is it needed? "Note that, since $\cP_n \subset \Z^n_+$, we are using the same notation $SST_{[n]}(\lambda)$ for $\lambda \in \Z^n_+$ as we use for $\lambda \in \cP_n $.")}
For example,\vskip 1mm
$$T=\resizebox{.18\hsize}{!}{${\def\lr#1{\multicolumn{1}{|@{\hspace{.6ex}}c@{\hspace{.6ex}}|}{\raisebox{-.3ex}{$#1$}}}\raisebox{-.6ex}
{$\begin{array}{ccc|ccc}
\cline{4-6}
& & & \lr{\, 1\,} & \lr{\, 1\,} & \lr{\, 2\,}  \\
\cline{4-6}
& &   & \lr{2}& \lr{3} & \\
\cline{4-5}
  &   &  &  &  &   \\
\cline{2-3}
  & \lr{ \hyphen 4 }  &  \lr{ \hyphen 3 }&  &  & \\
\cline{2-3}
\end{array}$}} $}\ \ \in \ \ \SST_{[4]}(3,2,0,-2)$$\vskip 2mm
\noindent with ${\rm wt}_{[4]}(T)=  2 \epsilon_1 + 2\epsilon_2 -\epsilon_4$.

For $0 \le t \le n$, let $T$ be a tableau in $\SST_{[n]}( 0^{n-t}, (-1)^t )$ with the entries $ -b_1, \ldots,  -b_t$. We denote
by $\sigma(T)$ the tableau in $\SST_{[n]}( 1^{n-t}, 0^t )$ with the entries $[n] \setminus \{b_1, \ldots, b_t \}$.
For an arbitrary tableau $T \in \SST_{[n]}(\lambda)$, by applying $\sigma$ to the $-1$st column of $T$, we have a bijection
\begin{equation} \label{Eq:sigma}
 \sigma: \SST_{[n]}(\lambda) \rightarrow \SST_{[n]}(\lambda + (1^n)),
\end{equation}
where ${\rm wt}_{[n]}(\sigma(T))={\rm wt}_{[n]}(T)+\sum_{i=1}^n\epsilon_i$.
Let $\lambda=(\lambda_1 ,\ldots, \lambda_n) \in \cP_n$  and  $T \in \SST_{[n]}(\lambda)$.
For $d \ge \lambda_1$, we set
$\delta_d(\lambda)  = (d^n) - (\lambda_n, \ldots, \lambda_1)$ and
$\delta_d (T) = (\sigma^{-d}(T))^\pi \in \SST_{[-n]^\pi}(\delta_d(\lambda))$.
Identifying $-k\in [-n]^{\pi}$ with $k\in [n]$,  we get a bijection
\begin{equation}\label{Eq:delta}
\delta_d: \SST_{[n]}(\lambda) \rightarrow \SST_{[n]}( \delta_d(\lambda) ).
\end{equation}

\begin{ex} \label{Ex: Ex2}
{\rm
Let $n=3$ and
$$
Q =
\resizebox{.115\hsize}{!}{${\def\lr#1{\multicolumn{1}{|@{\hspace{.6ex}}c@{\hspace{.6ex}}|}{\raisebox{-.3ex}{$#1$}}}\raisebox{-.0ex}
{$\begin{array}{cccc}
\cline{1-4}
 \lr{1} & \lr{1} & \lr{2} & \lr{2} \\
\cline{1-4}
\lr{2}  & \lr{3} &  \lr{3} &     \\
\cline{1-3}
\end{array}$}}$} \ \  \in \SST_{[3]}(4,3,0).
$$
Then we have $\delta^{-4}(4,3,0) = (4,1,0)$ and
$
\sigma^{-4}(Q) =
\resizebox{.125\hsize}{!}{${\def\lr#1{\multicolumn{1}{|@{\hspace{.6ex}}c@{\hspace{.6ex}}|}{\raisebox{-.3ex}{$#1$}}}\raisebox{-.0ex}
{$\begin{array}{cccc}
\cline{4-4}
 &  &  & \lr{\hyphen 3} \\
\cline{1-4}
\lr{\hyphen 3} & \lr{\hyphen 2} & \lr{\hyphen 1} & \lr{\hyphen 1} \\
\cline{1-4}
\end{array}$}}$}\ .
$
Thus
$$
\delta_4(Q) =
\resizebox{.125\hsize}{!}{${\def\lr#1{\multicolumn{1}{|@{\hspace{.6ex}}c@{\hspace{.6ex}}|}{\raisebox{-.3ex}{$#1$}}}\raisebox{-.0ex}
{$\begin{array}{cccc}
\cline{1-4}
\lr{1}  & \lr{1} & \lr{2} & \lr{3} \\
\cline{1-4}
\lr{3} &  \\
\cline{1-1}
\end{array}$}}$} \ \in \SST_{[3]}(4,1,0).
$$
}
\end{ex}

\subsection{Parabolically semistandard tableaux}\label{Sec:PSST}
Now, we review the notion of parabolically semistandard tableaux\footnote{These were called $\A/\B$-semistandard tableaux in \cite{K08}. }  introduced in \cite{K08} to study a combinatorial aspect of Howe dual pairs of type $A$.

Let $\lambda\in\mathbb{Z}_+^n$ be given. A {\it parabolically semistandard tableau of shape $\lambda$ with respect to $(\A,\B)$} is a pair of tableaux $(T^+,T^-)$ such that
\begin{equation*}
T^+\in \SST_{\A}((\lambda+(d^n))/\mu), \ \ \ \ T^- \in
\SST_{\B}((d^{n})/\mu),
\end{equation*}
for some integer $d\geq 0$ and $\mu\in\cP_n$ satisfying (1)
$\lambda+(d^n)\in\cP_n$, (2) $\mu\subset (d^n), \mu\subset \lambda+(d^n)$. We call $n$ the {\it level} of $T$ and define the weight of $T$ to be
\begin{equation*}
{\rm wt}_{\A/\B}(T)={\rm wt}_{\A}(T^+)-{\rm wt}_{\B}(T^-)\in P_{\A}\oplus P_{\B}.
\end{equation*}
We denote by $\SST_{\A/\B}(\lambda)$ the set of parabolically semistandard tableaux of shape $\lambda$.

Roughly speaking, $T\in \SST_{\A/\B}(\lambda)$ is a pair of an $\A$-semistandard tableau $T^+$ and a $\B$-semistandard tableau $T^-$, the difference of whose shapes is $\lambda$. For example, if $\A=\B=\Z_{>0}$, then the pair $(T^+,T^-)$ with
$$\resizebox{.65\hsize}{!}{$T^+\ \ =\ \ {\def\lr#1{\multicolumn{1}{|@{\hspace{.6ex}}c@{\hspace{.6ex}}|}{\raisebox{-.3ex}{$#1$}}}\raisebox{-.6ex}
{$\begin{array}{ccc|ccc}
& & & & & \\
\cline{3-6}
&   &  \lr{{\bf 1}} & \lr{1} & \lr{2} & \lr{2} \\
\cline{2-6}
  & \lr{\bf{1}}  & \lr{\bf{2}} &  \lr{2} & \lr{4} &   \\
\cline{1-5}
 \lr{\bf{2}} &  \lr{\bf{3}} &  \lr{\bf{3}} &  &  &   \\
\cline{1-3}
\lr{\bf{4}} & \cdot & \cdot &  &  &   \\
\cline{1-1}
& & & & &
\end{array}$}}
\ \ \ \ \ \ \ \ \ \  \ T^- \ \ = \ \
{\def\lr#1{\multicolumn{1}{|@{\hspace{.6ex}}c@{\hspace{.6ex}}|}{\raisebox{-.3ex}{$#1$}}}\raisebox{-.6ex}
{$\begin{array}{ccc|ccc}
& & & & & \\
\cline{3-3}
&   &  \lr{\bf{1}} & \cdot & \cdot & \cdot  \\
\cline{2-3}
  & \lr{\bf{1}}  & \lr{\bf{2}} & \cdot & \cdot &   \\
\cline{1-3}
 \lr{\bf{2}} &  \lr{\bf{2}} &  \lr{\bf{4}} &  &  &   \\
\cline{1-3}
\lr{\bf{3}} & \lr{3} & \lr{5} &  &  &   \\
\cline{1-3}
& & & & &
\end{array}$}}
$}$$
\noindent  belongs to $\SST_{\A/\B}((3,2,0,-2))$, where the vertical lines in $T^+$ and $T^-$ correspond to the one in the generalized partition $\lambda=(3,2,0,-2)$,
and the bold-faced entries  denote ones in the overlapping parts of ${\rm sh}(T^+)$ and ${\rm sh}(T^-)$. In this case, we have ${\rm
sh}(T^+)=\lambda+(3^4)/(2,1,0,0)$, and ${\rm
sh}(T^-)=(3^4)/(2,1,0,0)$. \vskip 2mm

Let us describe an analogue of RSK correspondence for parabolically semistandard tableaux \cite{K08}. From now on, we assume that $\A$ and $\B$ are disjoint sets.
Let
$$ \mathcal{F}_{\A/\B} = \bigsqcup_{k\in \Z} \SST_{\A/\B}(k) $$
be the set of all parabolically semistandard tableaux of level 1, and $\mathcal{F}_{\A/\B}^n$ its $n$-fold product.
Let $\mathbf{T} = ( T_1,\ldots, T_n ) \in \mathcal{F}_{\A/\B}^n$ be given with $T_i=(T^+_i, T^-_i)$. We associate a pair  $(P_{\bf T}, Q_{\bf T})$, where $P_{\bf{T}}$ is  a parabolically semistandard tableau  of level $n$ and $Q_{\bf T}$ is a rational semistandard tableau of rank $n$ determined by the following steps:

\begin{itemize}
\item[($\kappa$-1)]  Let
\begin{align*} \label{Eq: rhorow1}
 (P,Q)=\rhocol( (T_1^-)^\pi, \ldots, (T_n^-)^\pi ).
\end{align*}
%Recall that $\rhocol$ is given in \eqref{eq:rho} and $\pi$ denotes the rotation.
Put $T^- = P^\pi$ and write $\sh(T^-) = (d^n) / \mu$ for some $d \ge 0$ and $\mu \in \cP_n$.

\item[($\kappa$-2)] Let $Q^\vee =  \delta_d (Q)$, which is of shape $\mu$, and let $\nu=(\nu_1, \ldots, \nu_n)$, where $\mathrm{wt}_{[n]}(Q^\vee) = \sum_{i\in [n]} \nu_i \epsilon_i$. By \eqref{Eq:rho}, there exist unique $S_i \in \SST_{[n]}(\nu_i)$ for $1\leq i\leq n$ such that
$$ \rhorow(S_1, \ldots, S_n) = (H^\mu, Q^\vee) \in \SST_{[n]}(\mu) \times \SST_{[n]}(\mu)_\nu ,$$
where $H^\mu$ is the tableau of shape $\mu$ with weight $\sum_{i\in [n]}\mu_i\epsilon_i$.

\item[($\kappa$-3)] For $1\leq i\leq n$, put $U_i = S_i * T_i^+$, which is an $[n]*\A$-semistandard tableau. Using \eqref{Eq:rho} once again, we let
$$(U, U_R)= \rhorow(U_1, \ldots, U_n)  \in \SST_{[n]*\A}(\lambda + (d)^n) \times \SST_{[n]}(\lambda + (d)^n) ,$$
for some $\lambda \in \Z_+^n$.

\item[($\kappa$-4)] Since $i < a$ for all $i\in [n]$ and $a\in \A$ in $[n]*\A$, we have
$ U = H^\mu * T^+$ for some $T^+ \in \SST_\A( \lambda+(d)^n / \mu  )$. Finally, we define
\begin{equation*}
\begin{split}
& P_{\mathbf{T}} = (T^+, T^-) \in \SST_{\A / \B} (\lambda),\\
& Q_{\mathbf{T}} = \sigma^{-d}(U_R) \in \SST_{[n]}(\lambda).
\end{split}
\end{equation*}
\end{itemize}

% {\color{red} Shall we give an example for $\mathbf{T} \mapsto  (P_{\mathbf{T}}, Q_{\mathbf{T}})$ in case of $\A=\Z_{>0}$ and $\B=\Z_{<0}$ ?}
\begin{ex} \label{Ex: RSK}
{\rm
Let $\A = \Z_{\ge0}'$ and $\B = \Z_{<0}'$. Note that $\A$ and $\B$ have only elements of odd degree. %We identify $ (\Z_{<0}')^\pi $ with $\Z_{>0}'$ and vice versa.
Consider
$$ \mathbf{T} = (T_1, T_2, T_3) \in \SST_{\A/\B}(3) \times \SST_{\A/\B}(1) \times \SST_{\A/\B}(0) \subset \F_{\A/\B}^3, $$
where
\begin{align*}
T_1 &= ( T_1^+, T_1^- )=
(\
\resizebox{.17\hsize}{!}{${\def\lr#1{\multicolumn{1}{|@{\hspace{.6ex}}c@{\hspace{.6ex}}|}{\raisebox{-.1ex}{$#1$}}}\raisebox{ .1ex}
{$\begin{array}{|ccccc}
\cline{1-5}
\lr{ 0'}  & \lr{1'} & \lr{ 3'} & \lr{4'} & \lr{5'}  \\
\cline{1-5}
\end{array}$}} $}
\ , \
\resizebox{.08\hsize}{!}{${\def\lr#1{\multicolumn{1}{|@{\hspace{.6ex}}c@{\hspace{.6ex}}|}{\raisebox{-.1ex}{$#1$}}}\raisebox{ .1ex}
{$\begin{array}{|cc}
\cline{1-2}
\lr{\hyphen 4'}  & \lr{\hyphen 3'} \\
\cline{1-2}
\end{array}$}} $}  \
), \\
T_2 &= ( T_2^+, T_2^- )=  (\
\resizebox{.14\hsize}{!}{${\def\lr#1{\multicolumn{1}{|@{\hspace{.6ex}}c@{\hspace{.6ex}}|}{\raisebox{-.1ex}{$#1$}}}\raisebox{ .1ex}
{$\begin{array}{|cccc}
\cline{1-4}
\lr{0'}  & \lr{ 2' } & \lr{ 6' } & \lr{ 7' } \\
\cline{1-4}
\end{array}$}} $}
\ , \
\resizebox{.125\hsize}{!}{${\def\lr#1{\multicolumn{1}{|@{\hspace{.6ex}}c@{\hspace{.6ex}}|}{\raisebox{-.1ex}{$#1$}}}\raisebox{ .1ex}
{$\begin{array}{|ccc}
\cline{1-3}
\lr{ \hyphen 4' }  & \lr{ \hyphen 2' } & \lr{ \hyphen 1'}  \\
\cline{1-3}
\end{array}$}} $}
\  ), \\
T_3 &= (T_3^+, T_3^-) = (\
\resizebox{.07\hsize}{!}{${\def\lr#1{\multicolumn{1}{|@{\hspace{.6ex}}c@{\hspace{.6ex}}|}{\raisebox{-.1ex}{$#1$}}}\raisebox{ .1ex}
{$\begin{array}{|cc}
\cline{1-2}
  \lr{ 4' } & \lr{ 5' } \\
\cline{1-2}
\end{array}$}} $}
\ , \
\resizebox{.085\hsize}{!}{${\def\lr#1{\multicolumn{1}{|@{\hspace{.6ex}}c@{\hspace{.6ex}}|}{\raisebox{-.1ex}{$#1$}}}\raisebox{ .1ex}
{$\begin{array}{|cc}
\cline{1-2}
  \lr{ \hyphen 2' } & \lr{ \hyphen 1' } \\
\cline{1-2}
\end{array}$}} $}
\  ).
\end{align*}\vskip 1mm

Then, we have
\begin{align*}
&\rhocol( (T_1^-)^\pi, (T_2^-)^\pi, (T_3^-)^\pi ) =
\rhocol(\
\resizebox{.08\hsize}{!}{${\def\lr#1{\multicolumn{1}{|@{\hspace{.6ex}}c@{\hspace{.6ex}}|}{\raisebox{-.1ex}{$#1$}}}\raisebox{ .1ex}
{$\begin{array}{|cc}
\cline{1-2}
  \lr{\hyphen 3' } & \lr{\hyphen 4' } \\
\cline{1-2}
\end{array}$}} $}
\ ,\
\resizebox{.12\hsize}{!}{${\def\lr#1{\multicolumn{1}{|@{\hspace{.6ex}}c@{\hspace{.6ex}}|}{\raisebox{-.1ex}{$#1$}}}\raisebox{ .1ex}
{$\begin{array}{|ccc}
\cline{1-3}
 \lr{\hyphen 1' } & \lr{\hyphen 2' } & \lr{\hyphen 4' } \\
\cline{1-3}
\end{array}$}} $}
\ , \
\resizebox{.08\hsize}{!}{${\def\lr#1{\multicolumn{1}{|@{\hspace{.6ex}}c@{\hspace{.6ex}}|}{\raisebox{-.1ex}{$#1$}}}\raisebox{ .1ex}
{$\begin{array}{|cc}
\cline{1-2}
  \lr{\hyphen 1' } & \lr{\hyphen 2'} \\
\cline{1-2}
\end{array}$}} $}\ )
\\
& = (P,Q)  \  = \left( \
\resizebox{.17\hsize}{!}{${\def\lr#1{\multicolumn{1}{|@{\hspace{.6ex}}c@{\hspace{.6ex}}|}{\raisebox{-.3ex}{$#1$}}}\raisebox{-.0ex}
{$\begin{array}{cccc}
\cline{1-4}
 \lr{\hyphen 1'} & \lr{\hyphen 2'} & \lr{\hyphen 3'} & \lr{\hyphen 4'} \\
\cline{1-4}
\lr{\hyphen 1'}  & \lr{\hyphen 2'} &  \lr{\hyphen 4'} &     \\
\cline{1-3}
\end{array}$}}$}
\ , \
\resizebox{.115\hsize}{!}{${\def\lr#1{\multicolumn{1}{|@{\hspace{.6ex}}c@{\hspace{.6ex}}|}{\raisebox{-.3ex}{$#1$}}}\raisebox{-.0ex}
{$\begin{array}{cccc}
\cline{1-4}
 \lr{1} & \lr{1} & \lr{2} & \lr{2} \\
\cline{1-4}
\lr{2}  & \lr{3} &  \lr{3} &     \\
\cline{1-3}
\end{array}$}}$}
\ \right)
\in \SST_{(\Z_{<0}')^\pi}(4,3,0) \times \SST_{[3]}(4,3,0),
\end{align*}
which yield
$T^-=
\resizebox{.165\hsize}{!}{${\def\lr#1{\multicolumn{1}{|@{\hspace{.6ex}}c@{\hspace{.6ex}}|}{\raisebox{-.3ex}{$#1$}}}\raisebox{-.0ex}
{$\begin{array}{cccc}
\cline{2-4}
   & \lr{\hyphen 4'} & \lr{\hyphen 2'} & \lr{\hyphen 1'} \\
\cline{1-4}
\lr{\hyphen 4'} & \lr{\hyphen 3'}  & \lr{\hyphen 2'} &  \lr{\hyphen 1'}     \\
\cline{1-4}
\end{array}$}}$}
$ (see Example \ref{Ex: Ex1}). We choose $d = 4$ and $\mu=(4,1,0)$. So we have
$$
Q^\vee =
\resizebox{.125\hsize}{!}{${\def\lr#1{\multicolumn{1}{|@{\hspace{.6ex}}c@{\hspace{.6ex}}|}{\raisebox{-.3ex}{$#1$}}}\raisebox{-.0ex}
{$\begin{array}{cccc}
\cline{1-4}
\lr{1}  & \lr{1} & \lr{2} & \lr{3} \\
\cline{1-4}
\lr{3} &  \\
\cline{1-1}
\end{array}$}}$} \quad  \in \SST_{[3]}(\mu)
$$
with ${\rm wt}_{[3]}(Q^\vee) = 2 \epsilon_1 + \epsilon_2 + 2\epsilon_3$, which implies $\nu = (2,1,2)$\ (see Example \ref{Ex: Ex2}).
It follows from
$$
\rhorow(\
\resizebox{.06\hsize}{!}{${\def\lr#1{\multicolumn{1}{|@{\hspace{.6ex}}c@{\hspace{.6ex}}|}{\raisebox{-.1ex}{$#1$}}}\raisebox{ .1ex}
{$\begin{array}{|cc}
\cline{1-2}
\lr{ 1 }  & \lr{ 1 } \\
\cline{1-2}
\end{array}$}} $}
\ ,\
\resizebox{.03\hsize}{!}{${\def\lr#1{\multicolumn{1}{|@{\hspace{.6ex}}c@{\hspace{.6ex}}|}{\raisebox{-.1ex}{$#1$}}}\raisebox{ .1ex}
{$\begin{array}{|c}
\cline{1-1}
\lr{ 2 }   \\
\cline{1-1}
\end{array}$}} $}
\ , \
\resizebox{.06\hsize}{!}{${\def\lr#1{\multicolumn{1}{|@{\hspace{.6ex}}c@{\hspace{.6ex}}|}{\raisebox{-.1ex}{$#1$}}}\raisebox{ .1ex}
{$\begin{array}{|cc}
\cline{1-2}
\lr{ 1 }  & \lr{ 1 } \\
\cline{1-2}
\end{array}$}} $}
\ ) =
\left( \
\resizebox{.125\hsize}{!}{${\def\lr#1{\multicolumn{1}{|@{\hspace{.6ex}}c@{\hspace{.6ex}}|}{\raisebox{-.3ex}{$#1$}}}\raisebox{-.1ex}
{$\begin{array}{cccc}
\cline{1-4}
\lr{1}  & \lr{1} & \lr{1} & \lr{1} \\
\cline{1-4}
\lr{2} &  \\
\cline{1-1}
\end{array}$}}$}
\ , \
\resizebox{.125\hsize}{!}{${\def\lr#1{\multicolumn{1}{|@{\hspace{.6ex}}c@{\hspace{.6ex}}|}{\raisebox{-.3ex}{$#1$}}}\raisebox{-.1ex}
{$\begin{array}{cccc}
\cline{1-4}
\lr{1}  & \lr{1} & \lr{2} & \lr{3} \\
\cline{1-4}
\lr{3} &  \\
\cline{1-1}
\end{array}$}}$}
\ \right)
\  \in \SST_{[3]}(\mu) \times \SST_{[3]}(\mu)_\nu
$$
that
\begin{align*}
U_1 &=
\resizebox{.25\hsize}{!}{${\def\lr#1{\multicolumn{1}{|@{\hspace{.6ex}}c@{\hspace{.6ex}}|}{\raisebox{-.3ex}{$#1$}}}\raisebox{-.1ex}
{$\begin{array}{ccccccc}
\cline{1-7}
\lr{\,1\,}  & \lr{\,1\,} & \lr{0'} & \lr{1'} & \lr{3'} & \lr{4'} & \lr{5'} \\
\cline{1-7}
\end{array}$}}$}
\ , \\
U_2 &=
\resizebox{.18\hsize}{!}{${\def\lr#1{\multicolumn{1}{|@{\hspace{.6ex}}c@{\hspace{.6ex}}|}{\raisebox{-.3ex}{$#1$}}}\raisebox{-.1ex}
{$\begin{array}{ccccc}
\cline{1-5}
\lr{\,2\,}  &  \lr{0'} & \lr{2'} & \lr{6'} & \lr{7'} \\
\cline{1-5}
\end{array}$}}$}
\ , \\
U_3 &=
\resizebox{.145\hsize}{!}{${\def\lr#1{\multicolumn{1}{|@{\hspace{.6ex}}c@{\hspace{.6ex}}|}{\raisebox{-.3ex}{$#1$}}}\raisebox{-.1ex}
{$\begin{array}{cccc}
\cline{1-4}
\lr{\,1\,}  &  \lr{\, 1\,} & \lr{4'} & \lr{5'}  \\
\cline{1-4}
\end{array}$}}$}
\ .
\end{align*}\vskip 1mm

\noindent
Thus, we have $\rhorow(U_1, U_2, U_3) = (U, U_R)$, where \vskip 1mm
\begin{align*}
U =
\resizebox{.32\hsize}{!}{${\def\lr#1{\multicolumn{1}{|@{\hspace{.6ex}}c@{\hspace{.6ex}}|}{\raisebox{-.3ex}{$#1$}}}\raisebox{-.1ex}
{$\begin{array}{cccccccccc}
\cline{1-9}
\lr{\,1\,}  & \lr{\,1\,}  & \lr{\,1\,}  & \lr{\,1\,}  &  \lr{2'} & \lr{4'} & \lr{5'} & \lr{6'}  & \lr{7'} \\
\cline{1-9}
 \lr{\,2\,}  &  \lr{0'} & \lr{3'} & \lr{4'} & \lr{5'}   \\
\cline{1-5}
\lr{0'} & \lr{1'} \\
\cline{1-2}
\end{array}$}}$}
\ , \ \
U_R =
\resizebox{.32\hsize}{!}{${\def\lr#1{\multicolumn{1}{|@{\hspace{.6ex}}c@{\hspace{.6ex}}|}{\raisebox{-.3ex}{$#1$}}}\raisebox{-.1ex}
{$\begin{array}{cccccccccc}
\cline{1-9}
\lr{\,1\,}  & \lr{\,1\,}  & \lr{\,1\,}  & \lr{\,1\,}  & \lr{\,1\,}  & \lr{\,1\,}  & \lr{\,1\,}  & \lr{\,2\,} & \lr{\,2 \,} \\
\cline{1-9}
 \lr{\,2\,}  &  \lr{\,2\,}  & \lr{\,2\,}  & \lr{\,3\,}  & \lr{\,3\,}  \\
\cline{1-5}
\lr{\,3\,} & \lr{\,3\,} \\
\cline{1-2}
\end{array}$}}$}
\ .
\end{align*}
Since $d=4$ and $U = H^{\mu} * T^+$ (see Example \ref{Ex: Ex1}), we have $\lambda = (5,1,-2)$ and \vskip 1mm
\begin{align*}
P_{\mathbf{T}} &= \left(\
\resizebox{.32\hsize}{!}{${\def\lr#1{\multicolumn{1}{|@{\hspace{.6ex}}c@{\hspace{.6ex}}|}{\raisebox{-.3ex}{$#1$}}}\raisebox{-.1ex}
{$\begin{array}{cccc|cccccc}
& & & & & & & & \\
\cline{5-9}
   &    &    &    &  \lr{2'} & \lr{4'} & \lr{5'} & \lr{6'}  & \lr{7'} \\
\cline{2-9}
   &  \lr{\bf 0'} & \lr{\bf 3'} & \lr{\bf 4'} & \lr{ 5'}   \\
\cline{1-5}
\lr{\bf 0'} & \lr{\bf 1'} & & & & & & & \\
\cline{1-2}
& & & & & & & & \\
\end{array}$}}$} \ \ , \ \
\resizebox{.165\hsize}{!}{${\def\lr#1{\multicolumn{1}{|@{\hspace{.6ex}}c@{\hspace{.6ex}}|}{\raisebox{-.3ex}{$#1$}}}\raisebox{-.0ex}
{$\begin{array}{cccc|}
& & & \\
& & & \\
\cline{2-4}
   & \lr{\bf \hyphen 4'} & \lr{\bf \hyphen 2'} & \lr{\bf \hyphen 1'} \\
\cline{1-4}
\lr{\bf \hyphen 4'} & \lr{\bf \hyphen 3'}  & \lr{\hyphen 2'} &  \lr{\hyphen 1'}     \\
\cline{1-4}
& & & \\
\end{array}$}}$}
\ \ \ \ \ \ \right)  , \\
Q_{\mathbf{T}}  &= \ \
\resizebox{.25\hsize}{!}{${\def\lr#1{\multicolumn{1}{|@{\hspace{.6ex}}c@{\hspace{.6ex}}|}{\raisebox{-.3ex}{$#1$}}}\raisebox{-.1ex}
{$\begin{array}{cc|cccccc}
\cline{3-7}
  &  & \lr{\,1\,}  & \lr{\,1\,}  & \lr{\,1\,}  & \lr{\,2\,} & \lr{\,2 \,} \\
\cline{3-7}
   &    & \lr{\,3\,}   \\
\cline{1-3}
\lr{\hyphen 3 } & \lr{\hyphen 2} \\
\cline{1-2}
\end{array}$}}$}
  \ .
\end{align*}\vskip 2mm
 }
\end{ex}

\begin{thm} \cite[Theorem 4.1]{K08} \label{Thm: RSK}
The map $\mathbf{T} \mapsto  (P_{\mathbf{T}}, Q_{\mathbf{T}})$ gives a bijection
$$ \kappa_{\A/\B} :  \mathcal{F}_{\A/\B}^n \longrightarrow  \bigsqcup_{\lambda \in \cP_{\A/\B, n}} \SST_{\A/\B}(\lambda) \times \SST_{[n]}(\lambda),$$
where $\cP_{ {\A/\B}, n} = \{ \lambda \in \Z_+^n \mid \SST_{\A/ \B}(\lambda) \ne \emptyset\} $.
\end{thm}

For ${\bf T} = (T_1, \ldots, T_n)\in \F^n_{\A/\B}$, we put
$$
{\rm wt}_{\A/\B}({\bf T})=\sum_{i=1}^n {\rm wt}_{\A/\B}( T_i) \in P_{\A} \oplus P_{\B} , \quad
{\rm wt}_{[n]}({\bf T})= \sum_{i=1}^n m_i\epsilon_i \in P_{[n]},
$$
where $T_i\in \SST_{\A/\B}(m_i)$ for $1\leq i\leq n$. One can observe that
$$
{\rm wt}_{\A/\B}({\bf T}) = {\rm wt}_{\A/\B}(P_{\mathbf{T}}),\qquad
{\rm wt}_{[n]}({\bf T}) =  {\rm wt}_{[n]}(Q_{\mathbf{T}}),
$$
and hence $\kappa_{\A/\B}$ preserves the weights.

We define the character of $\SST_{\A/\B}(\lambda)$  to be
\begin{equation}\label{Eq:PSST char}
S^{\A/\B}_{\lambda}=\sum_{T\in \SST_{\A/\B}(\lambda)}{\bf
x}_{\A/\B}^{T},
\end{equation}
where ${\bf x}_{\A/\B}^T ={\bf x}_{\A}^{T^+}({\bf x}_{\B}^{T^-})^{-1}$ for $T=(T^+,T^-)\in \SST_{\A/\B}(\lambda)$. Then Theorem \ref{Thm: RSK} establishes the following Cauchy-type identity:
\begin{equation}\label{Eq:Cauchy}
\prod_{i\in [n]}\frac{\prod_{a\in \A_1}(1+x_a x_i)\prod_{b\in
\B_1}(1+ x_b^{-1} x_i^{-1})}{\prod_{a\in \A_0}(1- x_a x_i)\prod_{b\in
\B_0}(1- x_b^{-1} x_i^{-1})}=
\sum_{\lambda \in \cP_{\A/\B, n}}S^{\A/\B}_{\lambda}s_{\lambda}({\bf
x}_{[n]}).
\end{equation}
Here $s_{\lambda}({\bf x}_{[n]})$ is the Laurent Schur polynomial corresponding to $\lambda\in \Z_+^n$.

Note that, when $\B=\emptyset$, we have $\cP_{\A/\B,n}=\cP_{\A}\cap \cP_{n}$, and $S^{\A/\B}_\lambda=s_{\lambda}({\bf x}_{\A})$,  which is the usual (super) Schur function or polynomial corresponding to $\lambda$, and the identity \eqref{Eq:Cauchy} recovers the well-known Cauchy identity. So a non-trivial generalization of Schur functions or more interesting cases occur when both $\A$ and $\B$ are non-empty.

\subsection{Howe duality and irreducible characters}\label{Sec:IrrChar}
The notion of parabolically semistandard tableaux and its RSK with rational semistandard tableaux for the general linear Lie algebra $\gl_n$ gives a unified combinatorial interpretation of various dualities of $(\g,\gl_n)$, where $\g$ is a general linear Lie superalgebra associated to $(\A,\B)$. We assume that the base field is $\mathbb{C}$.

Let us explain it in more detail. For an arbitrary countable $\Z_2$-graded totally ordered set $S$, let $V_S$ be a superspace with basis $\{\,v_{s}\,|\,s\in S\,\}$, and let $\gl_S$ be the general linear Lie superalgebra spanned by the elementary matrices $E_{s s'}$ for $s, s'\in S$, where the parity of $E_{ss'}$ is given by the sum of the parities of $s$ and $s'$ (cf.~\cite{Kac}).

Now we consider $\g=\gl_{\mc C}$ with $\mc{C}=\B\ast\A$. Let $$\mathscr{F}=S(V_\A \oplus V_\B^\vee)$$ be the super symmetric algebra generated by $V_\A \oplus V_\B^\vee$, where $V_\B^\vee$ is the restricted dual space of $V_\B$.
Recall that $\mathscr{F}$ can be viewed as an irreducible module over a Clifford-Weyl algebra.
Following the arguments in \cite[Sections 5.1 and 5.4]{CW12} (cf. \cite{Fr,KacR2}), one can define a semisimple action of ${\mf g}$ on $\mathscr{F}$, and a semisimple action of $\gl_n$ or $GL_n$ on $\mathscr{F}^{\otimes n}$  for $n\geq 1$ such that $\mathscr{F}^{\otimes n}$ decomposes into a finite-dimensional $\gl_n$-modules. Then the actions of $\mf{g}$ and $\gl_n$ commute with each other, and furthermore the image of $\mf{g}$ in ${\rm End}_{\mathbb{C}}(\mathscr{F}^{\otimes n})$ generates ${\rm End}_{\gl_n}(\mathscr{F}^{\otimes n})$. Therefore, we have the following multiplicity-free decomposition as a $(\mf{g},\gl_n)$-module,
\begin{equation}\label{Eq:FockspDecomp-weak}
\mathscr{F}^{\otimes n}\cong \bigoplus_{\lambda \in H_{\g,n}} L_\g(\lambda)\otimes L_n(\lambda),
\end{equation}
for a subset $H_{\g,n}$ of $\Z_+^n$, where $L_n(\lambda)$ is an irreducible $\gl_n$-module with highest weight $\lambda \in H_{\g,n}$, and $L_\g(\lambda)$ is an irreducible ${\mf g}$-module corresponding to $L_n(\lambda)$.
We define the character ${\rm ch}L_\g(\lambda)$ to be the trace of the operator $\prod_{c\in \mc{C}}x_c^{E_{cc}}$ on $L_\g(\lambda)$ for $\lambda\in H_n$.

Then we have the following decomposition, which is often referred to as Howe duality (for type $A$) (cf. \cite{CL,CLZ,CW03,Fr,H,KacR2,KV}).

\begin{thm}\label{Thm:HoweDuality} Let $\A$ and $\B$ be given. For $n\geq 1$, we have
\begin{equation*}\label{Eq:FockspDecomp}
\mathscr{F}^{\otimes n}\cong \bigoplus_{\lambda \in \cP_{\A/\B,n}} L_\g(\lambda)\otimes L_n(\lambda),
\end{equation*}
as a $(\mf{g},\gl_n)$-module, that is, $H_{\g,n}=\cP_{\A/\B,n}$, and the irreducible character ${\rm ch}L_\g(\lambda)$ is given by $S_\lambda^{\A/\B}$ for $\lambda\in \cP_{\A/\B,n}$.
\end{thm}
\pf Consider the operator $D=\prod_{c\in \mc{C}}x_c^{E_{cc}}\prod_{i\in [n]}x_i^{e_{ii}}$, where $e_{ii}$ is the $i$-th elementary diagonal matrix in $\gl_n$. Taking the trace of $D$ on both sides of \eqref{Eq:FockspDecomp-weak}, we have
\begin{equation*}
\prod_{i\in [n]}\frac{\prod_{a\in \A_1}(1+x_a x_i)\prod_{b\in
\B_1}(1+ x_b^{-1} x_i^{-1})}{\prod_{a\in \A_0}(1- x_a x_i)\prod_{b\in
\B_0}(1- x_b^{-1} x_i^{-1})}=
\sum_{\lambda \in H_{\g,n}}{\rm ch}L_\g(\lambda)s_{\lambda}({\bf x}_{[n]}).
\end{equation*}
Thus by the Cauchy-type identity \eqref{Eq:Cauchy} and the linear independence of Laurent Schur polynomials, we conclude that $H_{\g,n}=\cP_{\A/\B,n}$ and $S_\lambda^{\A/\B}={\rm ch}L_\g(\lambda)$ for $\lambda  \in \cP_{\A/\B, n}$.
\qed\vskip 2mm

Note that $L_\g(\lambda)$'s are mutually non-isomorphic irreducible $\mf{g}$-modules for $\lambda\in \bigcup_{n\geq 1}\cP_{\A/\B,n}$ and the tensor product $L_\g(\mu)\otimes L_\g(\nu)$ for $\mu\in \cP_{\A/\B,m}$ and $\nu\in \cP_{\A/\B,n}$ decomposes into a direct sum of $L_\g(\lambda)$'s for $\lambda\in \cP_{\A/\B,m+n}$ with finite multiplicity given by a Littlewood-Richardson number (see \cite[Theorem 4.7]{K08}).
Also, $L_\g(\lambda)$ is semisimple over a maximal Levi subalgebra ${\mf l}=\gl_{\A}\oplus \gl_{\B}$ of $\mf{g}$, and expanding $S^{\lambda}_{\A/\B}$ as a linear combination of $s_{\mu}({\bf x}_{\A})s_{\nu}({\bf x}^{-1}_{\B})$ for $\mu, \nu\in \cP$ (see \cite[Proposition 3.14]{K08}) gives a branching rule with respect to ${\mf l}$ or its associated maximal parabolic sub algebra.

Recall that when $\A$ is finite with $\A=\A_0$ or $\A_1$ and $\B=\emptyset$, the decomposition in Theorem \ref{Eq:FockspDecomp} is the classical $(\gl_\ell,\gl_n)$-Howe duality on symmetric algebra or exterior algebra generated by $\mathbb{C}^\ell\otimes \mathbb{C}^n$, where $\ell=|\A|$ (cf.~\cite{H}).
Below we list some of important examples where both $\A$ and $\B$ are non-empty, and ${\mf g}$ is a usual general linear Lie algebra (see \cite{K08} for more detailed exposition).

\begin{ex} \label{Ex: Howe}
{ \rm
(1) If $(\A,\B)=(\Z'_{\ge0},  \Z_{<0}'  )$, then $S_\lambda^{\A/\B}$ ($\lambda\in \Z_+^n$) is the character of an integrable highest weight module over the general linear Lie algebra $\gl_\infty$ with  highest weight of positive level $n$. The identity \eqref{Eq:Cauchy} corresponds to the $(\gl_\infty,\gl_n)$-duality on the level $n$ fermionic Fock space $\mathscr{F}^{\otimes n}$ \cite{Fr}.
In particular, $\SST_{\A/\B}(k)$ $(k\in \Z)$ can be identified with a linear basis of the level $1$ fermionic Fock space $\mathscr{F}$ of charge $k$, which is realized by $\cP$ \cite[Section 1]{JM83}, by mapping an element
$$T = (\ \begin{array}{|c|c|c|c|} \hline
          i'_s  & \cdots & i'_2 & i'_1 \\ \hline
        \end{array}\ , \
        \begin{array}{|c|c|c|c|} \hline
          \hyphen j'_1  & \hyphen j'_2 & \cdots & \hyphen j'_r   \\ \hline
        \end{array} \
        ),
$$ where $i_p\geq 0$, $j_q>0$ and $s-r=k$,
to the following Young diagram
$$\lambda_T =
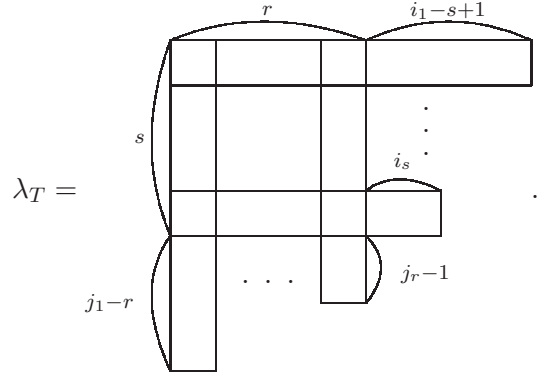
\xy
(0,20)*{};(48,20)*{} **\dir{-};
(0,14)*{};(48,14)*{} **\dir{-};
(0,0)*{};(36,0)*{} **\dir{-};
(0,-6)*{};(36,-6)*{} **\dir{-};
(0,-24)*{};(6,-24)*{} **\dir{-};
(20,-15)*{};(26,-15)*{} **\dir{-};
(0,20)*{};(0,-24)*{} **\dir{-};
(6,20)*{};(6,-24)*{} **\dir{-};
(20,20)*{};(20,-15)*{} **\dir{-};
(26,20)*{};(26,-15)*{} **\dir{-};
(36,0)*{};(36,-6)*{} **\dir{-};
(48,20)*{};(48,14)*{} **\dir{-};
(0,20)*{}; (0,-6)*{} **\crv{(-5,6)};
(0,-6)*{}; (0,-24)*{} **\crv{(-5,-13)};
(0,20)*{}; (26,20)*{} **\crv{(12,25)};
(26,20)*{}; (48,20)*{} **\crv{(37,25)};
(26,0)*{}; (36,0)*{} **\crv{(30,3)};
(26,-6)*{}; (26,-15)*{} **\crv{(30,-10)};
(-4,7)*{ _s}; (-8,-15)*{_{j_1-r}};
(13,24)*{ _r}; (37,24)*{_{i_1-s+1}};
(31,4)*{ _{i_s}}; (34,-11)*{ _{j_r - 1}};
(10,-12)*{ .}; (13,-12)*{ .};(16,-12)*{ .};
(34,5)*{ .}; (34,8)*{ .};(34,11)*{ .};
\endxy .
$$
For example, the Young diagrams corresponding to $T_1$, $T_2$ and $T_3$ in Example \ref{Ex: RSK} are
$\lambda_{T_1} = (3,3,3,2,2,2,2)$, $\lambda_{T_2} = (7,7,4,3,1)$ and $\lambda_{T_3} = (6,6)$, respectively.

(2) If $(\A,\B)=(\Z_{\ge 0},\Z_{<0})$, then $S_\lambda^{\A/\B}$ ($\lambda\in \Z_+^n$) is the character of an irreducible (non-integrable) highest weight module over $\gl_\infty$ with highest weight of negative level $-n$, which appears in the $(\gl_\infty,\gl_n)$-duality on the level $n$ bosonic Fock space \cite{KacR2}.

(3) If $(\A,\B)=([-q],[p])$ for $p,q\in \Z_{>0}$, then $S_\lambda^{\A/\B}$ ($\lambda\in \Z_+^n$) is equal to the character of an infinite-dimensional irreducible $\gl_{p+q}$-module, which is unitarizble. The family of irreducible representations appears in $(\gl_{p+q},\gl_n)$-duality on the symmetric algebra $S(\mathbb{C}^p\otimes \mathbb{C}^n \oplus \mathbb{C}^{q\ast}\otimes \mathbb{C}^{n\ast})$ \cite{KV}, which are called holomorphic discrete series or oscillator modules.
}
\end{ex}

\subsection{Hall-Littlewood functions}
Let $q$ be an indeterminate.
Fix $n\geq 1$. For $\mu\in \cP_n$, let
$P_\mu({\bf x}_{[n]},q)$
be the Hall-Littlewood polynomial in ${\bf x}_{[n]}$ associated to $\mu$ \cite[Chapter III.2]{Mac95}.
For $\mu\in \Z_+^n$, we put  $P_\mu({\bf x}_{[n]},q)=(x_1\ldots x_n)^{-d} P_{\mu+(d^n)}({\bf x}_{[n]},q)$ for some $d\geq 0$ such that $\mu+(d^n)\in \cP_n$,
which is independent of $d$ and hence well-defined.

Consider a formal power series $Q^{\A/\B}_{\lambda}$ in ${\bf x}_{\A}$ and ${\bf x}^{-1}_{\B}$, which is determined by the following Cauchy-type identity:
\begin{equation}\label{Eq:Cauchy-2}
\prod_{i\in [n]}\frac{\prod_{a\in \A_1}(1+x_a x_i)\prod_{b\in
\B_1}(1+ x_b^{-1} x_i^{-1})}{\prod_{a\in \A_0}(1- x_a x_i)\prod_{b\in
\B_0}(1- x_b^{-1} x_i^{-1})}=
\sum_{\lambda\in\Z_+^n}Q^{\A/\B}_{\lambda}P_{\lambda}({\bf
x}_{[n]},q).
\end{equation}

The following is a well-known identity:
\begin{equation}\label{Eq:s to P}
s_\lambda({\bf x}_{[n]})=\sum_{\mu\in \Z_+^n}K_{\lambda \mu}(q) P_\mu({\bf x}_{[n]},q),
\end{equation}
for $\lambda\in\Z_+^n$, where $K_{\lambda \mu}(q)$ ($\lambda, \mu \in \cP_n$) are the Kostka-Foulkes polynomials or Lusztig's $q$-weight multiplicities of type $A_{n-1}$. Here we set $K_{\lambda \mu}(q)=K_{\lambda+(d^n)\, \mu + (d^n)}(q)$ for $d\geq 1$ with $\lambda+(d^n), \mu + (d^n) \in \cP_n$, which is independent of the choice of $d$. By \eqref{Eq:Cauchy}, \eqref{Eq:Cauchy-2} and \eqref{Eq:s to P}, we have
\begin{equation}\label{Eq:Q function}
Q^{\A/\B}_\mu = \sum_{\lambda\in\Z_+^n}K_{\lambda \mu}(q)S^{\A/\B}_\lambda,
\end{equation}
for $\mu\in \Z_+^n$.
Since $K_{\lambda\mu}(q)$ has nonnegative integral coefficients with $K_{\lambda\mu}(1)=\left| \SST_{[n]}(\lambda)_\mu\right|$, we may view $Q^{\A/\B}_\mu $ as a $q$-analogue of the character of
\begin{equation}
\F^\mu_{\A/\B}= \SST_{\A/\B}(\mu_1)\times\cdots\times \SST_{\A/\B}(\mu_n),
\end{equation}
by Theorem \ref{Thm: RSK}. %\eqref{Eq:Cauchy}.
Recall that $Q^{\A/\B}_\mu $ is a modified Hall-Littlewood function for $\mu\in \cP$ when $\A=\Z_{>0}$ or $\Z_{>0}'$ and $\B=\emptyset$, and it can be realized as a graded character of a tensor product of KR crystals with respect to an energy function of affine type $A$ \cite{NaYa}.

Our main goal is to introduce a purely combinatorial statistic on $\F^\mu_{\A/\B}$, which realizes \eqref{Eq:Q function} as a graded character of $\F^{\mu}_{\A/\B}$ for arbitrary $\A$ and $\B$, also generalizing the usual energy functions on sequences of row (or column) tableaux.

\section{Affine crystals and charge statistic}
\subsection{Crystals} \label{Sec: Crystal}

Let us give a brief review on crystals (cf.\ \cite{HK,Kas94}).
Let $\g$ be the Kac-Moody algebra associated to a symmetrizable generalized Cartan matrix
$A =(a_{ij})_{i,j\in I}$. Let $P^\vee$ be the dual weight lattice,
$P = {\rm Hom}_\Z( P^\vee,\Z)$ the weight lattice,
$\Pi^\vee=\{\,h_i\,|\,i\in I\,\}$ the set of simple coroots, and
$\Pi=\{\,\alpha_i\,|\,i\in I\,\}$ the set of simple roots of $\g$ such that $\langle h_i,\alpha_j \rangle=a_{ij}$ for $i,j\in I$.
Let $U_q(\g)$ be the quantized enveloping algebra of $\g$.

A {\it $\g$-crystal} (or {\it crystal} for short) is a set
$B$ together with the maps ${\rm wt} : B \rightarrow P$,
$\varepsilon_i, \varphi_i: B \rightarrow \mathbb{Z}\cup\{-\infty\}$ and
$\te_i, \tf_i: B \rightarrow B\cup\{{\bf 0}\}$ ($i\in I$) satisfying certain axioms.
For a dominant integral weight $\Lambda$ for $\g$, we denote by $B(\Lambda)$ the crystal  associated to the irreducible highest weight $U_q(\g)$-module with highest weight
$\Lambda$.

For a crystal $B$, we denote its dual by $B^\vee$, which is a set $B^\vee = \{b^\vee \mid b \in B \}$ with
\begin{align*}
&{\rm wt}(b^\vee) = - {\rm wt}(b), \\
& \varepsilon_i(b^\vee) = \varphi_i(b), \quad \varphi_i(b^\vee) = \varepsilon_i(b),\\
& \te_i(b^\vee) = \tf_i(b)^\vee, \ \ \tf_i(b^\vee) = \te_i(b)^\vee,
\end{align*}
for $b\in B$ and $i\in I$.
A  tensor product $B_1\otimes B_2$  of crystals $B_1$ and $B_2$
is defined to be a crystal, which is $B_1\times B_2$  as a set with elements  denoted by
$b_1\otimes b_2$, satisfying  {\allowdisplaybreaks
\begin{equation*}
\begin{split}
{\rm wt}(b_1\otimes b_2)&={\rm wt}(b_1)+{\rm wt}(b_2), \\
\varepsilon_i(b_1\otimes b_2)&= {\rm
max}\{\varepsilon_i(b_1),\varepsilon_i(b_2)-\langle {\rm
wt}(b_1),h_i\rangle\}, \\
\varphi_i(b_1\otimes b_2)&= {\rm max}\{\varphi_i(b_1)+\langle {\rm
wt}(b_2),h_i\rangle,\varphi_i(b_2)\},\\
{\te}_i(b_1\otimes b_2)&=
\begin{cases}
{\te}_i b_1 \otimes b_2, & \text{if $\varphi_i(b_1)\geq \varepsilon_i(b_2)$}, \\
b_1\otimes {\te}_i b_2, & \text{if
$\varphi_i(b_1)<\varepsilon_i(b_2)$},
\end{cases}\\
{\tf}_i(b_1\otimes b_2)&=
\begin{cases}
{\tf}_i b_1 \otimes b_2, & \text{if  $\varphi_i(b_1)>\varepsilon_i(b_2)$}, \\
b_1\otimes {\tf}_i b_2, & \text{if $\varphi_i(b_1)\leq
\varepsilon_i(b_2)$},
\end{cases}
\end{split}
\end{equation*}
\noindent for $i\in I$. Here we assume that ${\bf 0}\otimes
b_2=b_1\otimes {\bf 0}={\bf 0}$.}

Given $b_1 \in B_1$ and $b_2 \in B_2$,
we write $b_1 \equiv b_2$ if there is an isomorphism of crystals $C(b_1) \rightarrow C(b_2)$ mapping $b_1$ to $b_2$, where $C(b_i)$ denotes the connected component of $b_i$ in $B_i$ for $i=1,2$.\vskip 2mm

\subsection{$A_{n-1}$-crystals}
Fix a positive integer $n \ge 2$.
Suppose that $\g = A_{n-1}$ or the associated generalized Cartan matrix is of type $A_{n-1}$ with $I=\{\,1,\ldots, n-1\,\}$. We assume that its weight lattice is $P_n:=P_{[n]}$. We often identify $\lambda \in \Z_+^{n}$ with the dominant integral weight $\sum_{i=1}^n\lambda_i\epsilon_i$. Let $\Delta_{n-1}^+=\{\,\epsilon_s-\epsilon_t\,|\,1\leq s<t\leq n\,\}$ the set of positive roots. The Weyl group is the symmetric group $\mf{S}_n$ on $n$ letters generated by the transposition $r_j=(j\ j+1)$ for $j=1,\ldots, n-1$.
From now on, we always denote the associated data of an $A_{n-1}$-crystal by $\tte_j$, $\ttf_j$, $\upvarepsilon_j$, $\upvarphi_j$ ($j=1,\ldots, n-1$) and ${\rm wt}$.

We may regard $[n]$ as $B(\epsilon_1)$ the crystal of the natural representation, and $[-n]$ as its dual. Given $\lambda\in \Z_+^n$, $\SST_{[n]}(\lambda)$ has an $A_{n-1}$-crystal structure by  regarding $w_{\rm col}(T)$ for $T\in \SST_{[n]}(\lambda)$ as an element in $[n]^{\otimes p}\otimes [-n]^{\otimes q}$ for some $p,q\geq 0$.
Here, we understand that $w_{\rm col}(T)$ is the word with letters in $[n]\cup [-n]$ obtained from $T$ by column reading as usual.
Then we have $\SST_{[n]}(\lambda)\cong B(\lambda)$ (cf. \cite{KN}).
It is not difficult to see that
\begin{align} \label{Eq:sigma delta}
\sigma^d (T) \equiv T,\qquad  \delta_d(S) \equiv S^{\vee},
\end{align}
for $T \in \SST_{[n]}(\lambda)$ with $\lambda\in \Z_+^n$ and $S \in \SST_{[n]}(\mu)$ with $\mu\in \cP_n$  up to a shift of weight by $d(\epsilon_1+\ldots+\epsilon_n)$ ($d \in \Z$) or as elements in $A_{n-1}$-crystals with the weight lattice $P_{n}/\Z(\epsilon_1+\ldots +\epsilon_n)$ (see \eqref{Eq:sigma} and \eqref{Eq:delta} for $\sigma$ and $\delta_d$).
The crystal equivalence $\equiv$ is also compatible with row and column insertions, that is,
$(T \leftarrow a) \equiv T \otimes a$ and $(a \rightarrow T ) \equiv a \otimes T$
for $a\in [n]$ and $T \in \SST_{[n]}(\lambda)$ with $\lambda\in \cP_n$.

\subsection{Charge statistic}
For $\lambda, \mu \in \cP_n$ and $T\in \SST_{[n]}(\lambda)_\mu$, we denote by $c(T)$ the {\it charge} of $T$, which was introduced by
Lascoux and Sch\"{u}tzenberger \cite{LS}.
It is shown in \cite{LS2} that
\begin{align} \label{Eq: KF poly}
K_{\lambda \mu}(q)=\sum_{T\in \SST_{[n]}(\lambda)_\mu}q^{c(T)}.
\end{align}

One can naturally induce a charge statistic on a regular $A_{n-1}$-crystal  $B$ as follows: Let $b\in B$ be given. First, note that the connected component $C(b)\subset B$ under $\tte_j$ and $\ttf_j$ ($j=1,\ldots,n-1$) is isomorphic to $B(\lambda)$ for some $\lambda\in \Z_n^+$. Choose $d\geq 0$ such that  $\lambda+(d^n)\in \cP_n$.  Since $B(\lambda)\cong \SST_{[n]}(\lambda+(d^n))$  as a $\{1,\ldots,n-1\}$-colored oriented graph by \eqref{Eq:sigma delta}, $b$ can be identified with a tableau $T\in \SST_{[n]}(\lambda+(d^n))$. Then we define
\begin{equation}
{\rm charge}(b)=c(T'),
\end{equation}
where $T'$ is a unique tableau with dominant weight in the $\mf{S}_n$-orbit of $T$. By definition of Lascoux and Sch\"{u}tzenberger's charge, it is not difficult to see that ${\rm charge}(b)$ does not depend on the choice of $d$.
In particular, we define for $\lambda,\mu\in \Z_+^n$
\begin{equation}
K_{\lambda \mu}(q)=\sum_{b\in B(\lambda),\ {\rm wt}(b)=\mu}q^{{\rm charge}(b)},
\end{equation}
which is equal to the usual Kostka-Foulkes polynomial $K_{\lambda+(d^n)\mu+(d^n)}(q)$ for some $d\geq 0$.

Recall that there is an intrinsic characterization of the charge statistic  \cite{LLT}, which is described only in terms of the geometry of the crystal graph $B(\lambda)$ for $\lambda\in \Z_+^n$. In Section \ref{Subsec:Howe}, we give another intrinsic characterization, which plays a crucial role in this paper. For this, we need the following statistic on a regular $A_{n-1}$-crystal $B$: for $\alpha=\epsilon_s-\epsilon_t\in \Delta_{n-1}^+$ and $b\in B$
\begin{equation} \label{Eq: epsilon and phi}
\begin{split}
\upvarepsilon_\alpha(b)&
=\upvarepsilon_s(S_{s+1}S_{s+2}\cdots S_{t-1}(b)),\\
\upvarphi_\alpha(b)&=\upvarphi_s(S_{s+1}S_{s+2}\cdots S_{t-1}(b)),
\end{split}
\end{equation}
where $S_j$ is the $\mf{S}_n$-action on $B$ associated to $r_j$.
Since $\upvarphi_\alpha(b)-\upvarepsilon_\alpha(b)=\langle {\rm wt}(b), \alpha^\vee \rangle $, where $\alpha^\vee=h_s+\cdots+h_{t-1}$ is the coroot of $\alpha$, one may think of $\upvarepsilon_\alpha(b)$ and $\upvarphi_\alpha(b)$ as information on {\it an} $\mathfrak{sl}_2$-string of $b$ with respect to $\alpha=r_{t-1}\ldots r_{s+1}(\alpha_s)$. We should remark that they depend on the choice of a simple root conjugate to $\alpha$. Here we choose it as $\alpha_s$.

\subsection{Affine $A_{\ell-1}^{(1)}$-crystals and energy function} \label{Sec: affine energy}
Fix a positive integer $\ell \ge 2$.
Suppose that $\g=A_{\ell-1}^{(1)}$ with $I=\{\,0,\ldots,\ell-1\,\}$  and
$\g_0=A_{\ell-1}$ is the subalgebra of $\g$ corresponding to $I \setminus \{ 0\}$.
For $1\leq r\leq \ell-1$, let $\varpi_r$ be the fundamental weight for $\g_0$ corresponding to the simple root $\alpha_r$. For $s\geq 1$, let $B^{r,s}$ denote the {\it Kirillov-Reshetikhin crystal} (or KR crystal for short) of type $A_{\ell-1}^{(1)}$, which is isomorphic to $B(s\varpi_r)$ as an $A_{\ell-1}$-crystal \cite{KMN2,S}. Let $u_{r,s}$ be the unique element in $B^{r,s}$ of weight $s\varpi_r$. For convenience, let us assume that $B^{0,s}$ and $B^{\ell,s}$ are trivial crystals.

Let $B_1$ and $B_2$ be two KR crystals with the classical highest weight elements $u_1$ and $u_2$, respectively. Let $\sigma=\sigma_{B_1,B_2} : B_1\otimes B_2 \longrightarrow B_2\otimes B_1$ be  a unique $A_{\ell-1}^{(1)}$-crystal isomorphism  called the {\it combinatorial $R$-matrix}. There exists a function $H=H_{B_1,B_2} : B_1\otimes B_2 \longrightarrow \Z$ such that $H$ is constant on each connected component in $B_1\otimes B_2$ as an $A_{\ell-1}$-crystal and
\begin{equation*}
\begin{split}
&H(\te_0(b_1\otimes b_2))=\\
&\begin{cases}
H(b_1\otimes b_2) +1, & \text{if $\te_0(b_1\otimes b_2)=\te_0(b_1)\otimes b_2$ and $\te_0(b'_2\otimes b'_1)=\te_0(b'_2)\otimes b'_1$} ,\\
H(b_1\otimes b_2) -1, &\text{if $\te_0(b_1\otimes b_2)=b_1 \otimes \te_0(b_2)$ and $\te_0(b'_2\otimes b'_1)= b'_2 \otimes \te_0(b'_1)$},\\
H(b_1\otimes b_2), & \text{otherwise},\\
\end{cases}
\end{split}
\end{equation*}
 for $b_1\otimes b_2\in B_1\otimes B_2$ with $b'_2\otimes b'_1=\sigma(b_1\otimes b_2)$. It is well-known that $H$ is unique up to an additive constant and is called the {\it local energy function} on $B_1\otimes B_2$ \cite{KMN2-2}. %\cmt{(Park: I cannot find the exact definition in this paper. The definition in Sec.\ 4.1 of the paper seems to be different a little. )}

Suppose that $B=B_1\otimes \cdots \otimes B_n$ is a tensor product of KR crystals. For $1\leq i\leq n-1$, let $\sigma_{i}$ be the $A_{\ell-1}^{(1)}$-crystal isomorphism of $B$, which acts as $\sigma_{B_i, B_{i+1}}$ on $B_i\otimes B_{i+1}$  and as identity elsewhere, and let $H_{i}$ be the function on $B$ given by $H_i(b_1\otimes \cdots \otimes b_n)=H_{B_i,B_{i+1}}(b_i\otimes b_{i+1})$.
The {\it energy function} $D_B : B \longrightarrow \Z$ is defined to be
\begin{equation}\label{Eq:Energy D}
D_B(b)= \sum_{1\leq i<j\leq n}H_{i}(\sigma_{i+1}\sigma_{i+2}\cdots \sigma_{j-1}(b)) \ \ \ \ \ (b\in B),
\end{equation}
which plays a very important role in the study of finite affine crystals (cf.~\cite{HKOTY,HKOTT}).
Note that $D_B$ is constant on each connected component in $B$ as an $A_{\ell-1}$-crystal, which therefore gives a natural $q$-analogue of the branching multiplicities
with respect to $A_{\ell-1}\subset A_{\ell-1}^{(1)}$.

%Note that $D_B$ is invariant under $\sigma_i$ for $1\leq i\leq n-1$, that is, $D_B=D_B\circ \sigma_i$ \cite{???}.

\subsection{Crystal skew Howe duality}\label{Subsec:Howe}

Let ${\bf M}_{\ell\times n}(\Z_2)$ be the set of $\ell\times n$ matrices ${\bf m}=(m_{ij})$ such that $m_{ij}=0,1$ for $1\leq i\leq \ell$ and $1\leq j\leq n$.

For $1\leq i\leq \ell$, let ${\bf m}_{(i)}$ denote the $i$th row of ${\bf m}$. We may identify each ${\bf m}_{(i)}$ with an $[n]$-semistandard tableau of single column whose entries are the column indices $j$ with $m_{ij}=1$, and hence regard ${\bf M}_{\ell\times n}(\Z_2)$ as an $A_{n-1}$-crystal by identifying ${\bf m}$ with ${\bf m}_{(\ell)}\otimes \cdots \otimes{\bf m}_{(1)}$ with respect to $\tte_j$, $\ttf_j$, $\upvarepsilon_j$, $\upvarphi_j$ for $1\leq j\leq  n-1$.

For $1\leq j\leq n$, let ${\bf m}^{(j)}$ denote the $j$th column of ${\bf m}$. In the same way, we  regard ${\bf M}_{\ell \times n}(\Z_2)$ as an $A_{\ell-1}$-crystal by identifying ${\bf m}$ with ${\bf m}^{(1)}\otimes \cdots \otimes{\bf m}^{(n)}$ with respect to $\te_i$, $\tf_i$, $\varepsilon_i$, $\varphi_i$ for $1\leq i\leq \ell-1$.

Then ${\bf M}_{\ell\times n}(\Z_2)$ is an $(A_{\ell-1},A_{n-1})$-bicrystal, that is, $\widetilde{x}_i \widetilde{\mathsf x}_j=\widetilde{\mathsf x}_j \widetilde{x}_i$ for all $i,j$, $x=e, f$ and $\mathsf{x} =\mathsf{e}, \mathsf{f}$, and the well-known (dual) RSK correspondence
\begin{equation}\label{Eq:RSK}
{\bf M}_{\ell\times n}(\Z_2) \longrightarrow \bigsqcup_{\lambda\in\cP} \SST_{[\ell]}(\lambda)\times \SST_{[n]}(\lambda')
\end{equation}
is a bicrystal isomorphism \cite{DK,K07}. Here $\lambda'$ denotes the conjugate of $\lambda$. This can be viewed as a crystal version of {\it skew $(\gl_\ell,\gl_n)$-Howe duality} (cf. \cite{H}).

Moreover, ${\bf M}_{\ell\times n}(\Z_2)$ is an $A_{\ell-1}^{(1)}$-crystal with respect to $\te_i$ and $\tf_i$ for $0\leq i\leq  \ell-1$, since each column of ${\bf m}\in {\bf M}_{\ell\times n}(\Z_2)$ can be considered as an element in a KR crystal $B^{r,1}$ for some $1\leq r\leq \ell-1$ or a trivial crystal. Note that $\te_0$ and $\tf_0$ do not commute with $\tte_j$ and $\ttf_j$ for $1\leq j\leq  n-1$, in general.

The affine $A_{\ell-1}^{(1)}$-crystal ${\bf M}_{\ell\times n}(\Z_2)$ is a union of tensor product of KR-crystals $B_1\otimes\cdots\otimes B_n$.
Thus, we can define the energy function $D$ on ${\bf M}_{\ell\times n}(\Z_2)$ as in \eqref{Eq:Energy D}, where we normalize the local energy function by requiring $H_{B_i, B_{i+1}}(u_i\otimes u_{i+1})=0$ for KR crystals $B_i$ with the classical highest weight elements $u_i\in B_i$ for $1\leq i\leq n$.
Then via \eqref{Eq:RSK} we can rewrite $D$ in terms of statistics on $A_{n-1}$-crystal as follows.

\begin{prop}\label{prop:energy description} For ${\bf m}\in {\bf M}_{\ell \times n}(\Z_2)$, we have
\begin{equation*}
D({\bf m})=-\sum_{\alpha\in \Delta^+_{n-1}}\min\{ \upvarepsilon_\alpha({\bf m}), \upvarphi_\alpha({\bf m})\}.
\end{equation*}
\end{prop}
\pf We may assume that ${\bf m}\in B=B_1\otimes \cdots \otimes B_n$, where $B_j = \SST_{[\ell]}(1^{t_j})=B^{t_j,1}$ for some $t_j$ ($1\leq j\leq  n$). We can check  in a straightforward manner by using the bicrystal structure on ${\bf M}_{\ell \times n}(\Z_2)$ that for $1\leq j\leq n-1$,
\begin{equation} \label{Eq: energy1}
\begin{aligned}
&\text{$\cdot$ $\sigma_{B_j,B_{j+1}}$ on $B_j\otimes B_{j+1}$ coincides with the $\mf{S}_n$-action $S_j$ on ${\bf m}$,}\\
&\text{$\cdot$ $H_{B_j,B_{j+1}}({\bf m}^{(j)}\otimes {\bf m}^{(j+1)})= - \min\{\,\upvarepsilon_j(  {\bf m}  ),\upvarphi_j( {\bf m} )\,\}  $,}
\end{aligned}
\end{equation}
(cf.~\cite[Section 3.5]{NaYa}), where
the second statement of $\eqref{Eq: energy1}$ may be understood as a crystal-theoretic interpretation of \cite[Rule 3.10]{NaYa}.
Then it follows from \eqref{Eq: epsilon and phi} and \eqref{Eq: energy1} that
\begin{equation*}
H_{s}(\sigma_{s+1}\sigma_{s+2}\cdots \sigma_{t-1}({\bf m}))= - \min\{ \upvarepsilon_\alpha({\bf m}), \upvarphi_\alpha({\bf m})\}  \ \ \ (1\leq s<t\leq n),
\end{equation*}
where $\alpha=\epsilon_s-\epsilon_t$. Hence, we get $D({\bf m})=-\sum_{\alpha\in \Delta^+_{n-1}}\min\{ \upvarepsilon_\alpha({\bf m}), \upvarphi_\alpha({\bf m})\}.$
\qed\vskip 2mm

Combining with the result of  Nakayashiki and Yamada  \cite{NaYa} (see also \cite{SW99,S} for its generalisation),
we obtain the following intrinsic characterization of charge statistic on a regular $A_{n-1}$-crystal.

\begin{thm} \label{Thm: charge}
Let $B$ be a regular $A_{n-1}$-crystal. For $b\in B$, we have
\begin{equation*}
{\rm charge}(b) =\sum_{\alpha\in \Delta^+_{n-1}}\min\{ \upvarepsilon_\alpha( b ), \upvarphi_\alpha(b)\}.
\end{equation*}
\end{thm}
\pf Given $b\in B$, we may assume that $b\in \SST_{[n]}(\lambda)$ for some $\lambda\in \cP_n$ up to a shift of its weight by $d(\epsilon_1+\ldots+\epsilon_n)$ ($d\geq 0$), say $b=T$. Let ${\bf m}$ be the unique matrix in ${\bf M}_{\ell \times n}(\Z_2)$ such that ${\bf m}_{(i)}$ corresponds to the $i$th column of
$T$ from the left-most column of $\lambda$. Since ${\rm charge}(b)$ is invariant under the Weyl group action, we may also assume that ${\rm wt}(b)$ is dominant, which corresponds to a partition  $\mu=(\mu_1,\ldots,\mu_n)$. Then ${\bf m}\in B^{\mu_1,1}\otimes \cdots \otimes  B^{\mu_n,1}$ as an $A_{\ell-1}^{(1)}$-crystal. Since we have $c(T)=- D({\bf m})$ by \cite[Section 4.1]{NaYa}, we have by Proposition \ref{prop:energy description}, $${\rm charge}(b)=c(T)=- D({\bf m})=\sum_{\alpha\in \Delta^+_{n-1}}\min\{ \upvarepsilon_\alpha(b), \upvarphi_\alpha(b)\}.$$ This completes the proof.
\qed

%\begin{rem}{\rm  One advantage of Theorem \ref{Thm: charge} is that the charge statistic does not depend on a realization of crystals. There is another intrinsic description of charge statistic in terms of $\mf{sl}_2$-strings,  which was introduced by Lascoux-Leclerc-Thibon \cite{LLT}.

% while our description is in terms of the $\mf{sl}_2$-strings over all positive roots. For $\lambda\in \cP_n$ and $T\in SST_{[n]}(\lambda)$, let
%\begin{equation*}
 %d(T)=\sum_{i=1}^{n-1}i \min\{\upvarepsilon_i(T), \upvarphi_i(T)\}.
%\end{equation*}
%Let $b(T)$ be the arithmetic mean of $d(T')$ for all $T'$ in the $\mf{S}_n$-orbit of $T$. Then it is shown by Lascoux-Leclerc-Thibon \cite{LLT} that
%\begin{equation}\label{Eq:LLT}
%b(T)=c(\eta (T)),
%\end{equation}
%where $\eta$ is the Lusztig's involution or Sch\"{u}tzenberger's involution on $SST_{[n]}(\lambda)$.
%}
%\end{rem}

\section{A combinatorial energy function}
In this section, we introduce a combinatorial energy function ${\mathsf D}$ on $\F^{\mu}_{\A/\B}$ for $\mu\in \Z_+^n$, which realizes $Q^\mu_{\A/\B}$ in \eqref{Eq:Q function} as a graded character of $\F^{\mu}_{\A/\B}$.

\subsection{Combinatorial energy function}
Consider
$\SST_{\A/\B}(k_1)\times \SST_{\A/\B}(k_2)$
for $k_1, k_2\in \Z$. Let $T_j= (T^+_j, T^-_j) \in \SST_{\A/\B}(k_j)$ be given for $j=1,2$ with
\begin{equation*}
{\rm wt}_{\A/\B}(T_j)=\sum_{a\in \A}m_{a j}\epsilon_a - \sum_{b\in \B}m_{bj}\epsilon_b.
\end{equation*}

First, we define a {\it local energy function}
$${\mathsf H} : \SST_{\A/\B}(k_1)\times \SST_{\A/\B}(k_2) \longrightarrow \Z,$$ following the steps below:

\begin{itemize}
\item[(${\mathsf H}$-1)]
Choose finite subsets $\A^\circ \subset \A$ and $\B^\circ \subset \B$ such that $T_j \in \SST_{\A^\circ/\B^\circ}(k_j)$ for $j=1,2$.
To each $i \in \A^\circ \sqcup \B^\circ$, we assign a sequence of $\pm$ signs as follows:

$$ \mathsf{s}_i = \left\{
                    \begin{array}{ll}
                      \underbrace{- \ \cdots\  -}_{m_{i2}}  \ \underbrace{+\ \cdots\ +}_{m_{i1}}  & \hbox{ if } i \in \A_0, \\
                      \underbrace{+}_{m_{i1}} \ \underbrace{-}_{m_{i2}}  & \hbox{ if } i \in \A_1, \\
                      \underbrace{- \ \cdots\  -}_{m_{i1}} \ \underbrace{+\ \cdots\ +}_{m_{i2}}  & \hbox{ if } i \in \B_0, \\
                      \underbrace{+ }_{m_{i2}}  \ \underbrace{- }_{m_{i1}}  & \hbox{ if } i \in \B_1.
                    \end{array}
                  \right.
$$

\item[(${\mathsf H}$-2)]  Let
$ \mathsf{s} = \mathsf{s}_{T_1,T_2}= (  \mathsf{s}_{a_k}\, \mathsf{s}_{a_{k-1}} \, \ldots \mathsf{s}_{a_1}\, \mathsf{s}_{b_1}\, \mathsf{s}_{b_2}\, \ldots\, \mathsf{s}_{b_l} )$
be their concatenation where $\A^\circ=\{ a_k > \cdots > a_1 \}$ and $\B^\circ=\{ b_1 < \cdots < b_l \}$, and cancel out all possible $(+\ -)$ pairs in $\mathsf{s}$ as far as possible to obtain a reduced sequence
$$
\mathsf{s}^{\rm red}= ( \underbrace{- \ \cdots\  -}_{\upvarepsilon}  \ \underbrace{+\ \cdots\ +}_{\upvarphi} ).
$$
Then we define
\begin{align*} \label{Eq: energy}
{\mathsf H} ( T_1, T_2 ) = - \min\{ \upvarepsilon, \upvarphi  \}.
\end{align*}
\end{itemize}

\begin{lem}\label{Lem:phi-epsilon}
With the same notations as above, we have $\upvarphi-\upvarepsilon = k_1-k_2$.
\end{lem}
\pf Let $p$ (resp. $q$) be the total number of $+$'s (resp. $-$'s) in ${\mathsf s}$. Then
$$p=\sum_{a\in \A}m_{a1} + \sum_{b\in \B}m_{b2}, \ \ \ q= \sum_{a\in \A}m_{a2}+\sum_{b\in \B}m_{b1}.$$ Since $\upvarphi - \upvarepsilon  = p-q$ and $k_j = \sum_{a\in \A}m_{a j}-\sum_{b\in \B}m_{b j}$ for $j=1,2$, we have
$\upvarphi-\upvarepsilon = k_1-k_2$.
\qed

Next, we define a {\it combinatorial $R$-matrix}
\begin{align*}
\upvarsigma : \SST_{\A/\B}(k_1) &\times \SST_{\A/\B}(k_2) \longrightarrow \SST_{\A/\B}(k_2) \times \SST_{\A/\B}(k_1)\ , \\
 (T_1&, T_2) \ \ \ \ \ \ \ \ \ \ \ \ \ \, \longmapsto \ \ \ \ \ \ \ \ \ \ \ \ \ (T'_2,T'_1)
\end{align*}
where $(T'_2,T'_1)$ is given by moving and rearranging some of the entries in $T_1$ and $T_2$ in the following way:
\begin{itemize}
\item[($\upvarsigma$-1)] If $k_1=k_2$, then put $(T'_2,T'_1)=(T_1,T_2)$.

\item[($\upvarsigma$-2)]
If $k_1 > k_2 $, then let  $y_1,\ldots,y_{k_1-k_2}$ be the entries in $T_1$ or $T_2$ corresponding to the first $k_1-k_2$ signs of $+$ in $\mathsf{s}^{\rm red}$ from the left (see Lemma \ref{Lem:phi-epsilon}).
For each $1\leq k\leq k_1-k_2$, if  $y_k\in\A$ (resp. $y_k\in\B$),  i.e.  $\boxed{y_k}$ appears $T_1^+$ (resp. $T_2^-$), then we move it to $T_2^+$ (resp. $T_1^-$) and rearrange the entries with respect to the total order on $\A$ (resp. $\B$).

\item[($\upvarsigma$-3)]
If $k_1 < k_2$, then let  $x_{k_2-k_1},\ldots,x_{1}$ be the entries in $T_1$ or $T_2$ corresponding to the first $k_2-k_1$ signs of $-$ in $\mathsf{s}^{\rm red}$ from the right.
For each $1\leq k\leq k_2-k_1$, if  $x_k\in\A$ (resp. $x_k\in\B$),  i.e. $\boxed{x_k}$ appears $T_2^+$ (resp. $T_1^-$), then we move it to $T_1^+$ (resp. $T_2^-$) and rearrange the entries with respect to the total order on $\A$ (resp. $\B$).
\end{itemize}

By definition, it is clear that $\upvarsigma\circ\upvarsigma = {\rm id}$.

\begin{ex} \label{Ex: energy}
{\rm
(1) If $\A  = [\ell]'$ and $\B = \emptyset$, then we have $\SST_{\A/\B}(k_j) = B^{k_j,1}$ for $j=1,2$, and identify $\SST_{\A/\B}(k_1) \times \SST_{\A/\B}(k_2)$ with $B^{k_1,1} \otimes B^{k_2,1}$.
As we have seen in Section \ref{Subsec:Howe}, we regard ${\bf M}_{\ell\times 2}(\Z_2)$ as a union of tensor products $B^{r,1} \otimes B^{s,1}$.
Since ${\bf M}_{\ell\times 2}(\Z_2)$ is an $(A_{\ell-1}, A_1)$-bicrystal,
we can apply $\upvarepsilon_1$ and $\upvarphi_1$ to $(T_1,T_2)$ in $\SST_{\A/\B}(k_1) \times \SST_{\A/\B}(k_2)$ (as an $A_1$-crystal). It follows from the definition of ${\mathsf H}$ that
$$ {\mathsf H}(T_1, T_2) = - \min \{ \upvarepsilon_1(T_1, T_2), \upvarphi_1(T_1, T_2)  \}.$$
Then we have  from \eqref{Eq: energy1} that
$${\mathsf H}(T_1, T_2) = H_{B^{k_1,1},B^{k_2,1}} (T_1 \otimes T_2).$$
Hence ${\mathsf H}$ on $B^{k_1,1}\times B^{k_2,1}$ coincides with the local energy function $H$ on $B^{k_1,1}\otimes B^{k_2,1}$ normalized by $H(u_{k_1,1}\otimes u_{k_2,1})=0$.
%\cmt{(Kwon: Revised some sentences slightly)}

%(1) If $\A  = [\ell]'$ and $\B = \emptyset$, then we have $SST_{\A/\B}(k_j) = B^{k_j,1}$ for $j=1,2$, and ${\mathsf H}$ on $B^{k_1,1}\times B^{k_2,1}$ \cmt{coincides with} the local energy function $H$ on $B^{k_1,1}\otimes %B^{k_2,1}$ normalized by $H(u_{k_1,1}\otimes u_{k_2,1})=0$ (see Proposition \ref{prop:energy description}).

(2)
If $\A  = [\ell]$ and $\B = \emptyset$, then we have $\SST_{\A/\B}(k_j) = B^{1, k_j}$ for $j=1,2$, and identify $\SST_{\A/\B}(k_1) \times \SST_{\A/\B}(k_2)$ with $B^{1,k_1} \otimes B^{1,k_2}$. As we have seen (1),
we consider the crystal version of ($\mathfrak{gl}_n, \mathfrak{gl}_2$)-Howe duality on the union of tensor products $B^{1,r} \otimes B^{1,s}$. Then, by the definition of ${\mathsf H}$, one can show that
$$ {\mathsf H}(T_1, T_2) = - \min \{ \upvarepsilon_1(T_1, T_2), \upvarphi_1(T_1, T_2)  \}, $$
for $(T_1, T_2) \in \SST_{\A/\B}(k_1) \times \SST_{\A/\B}(k_2)$. Interpreting \cite[Rule 3.11]{NaYa} from a point of view of crystal bases theory as $\eqref{Eq: energy1}$,
we can check that the map $-{\mathsf H}$ on $\SST_{\A/\B}(k_1) \times \SST_{\A/\B}(k_2)$ coincides with the local energy function $H$ on  $B^{1, k_2}\otimes B^{1, k_1}$ (in reverse order) normalized by $H(u_{1,k_2}\otimes u_{1,k_1})=\min\{k_1,k_2\}$ (cf.\ \cite[Section 3]{NaYa}). In both cases (1) and (2), $\upvarsigma$ is equal to the combinatorial $R$-matrix $\sigma$.

%\cmt{(Park: Maybe I am confusing... but I feel that
%\begin{align*}
%\text{${\mathsf H}$ on $B^{1, k_1}\otimes B^{1, k_2}$ = $-H$ on  $B^{1, k_2}\otimes B^{1, k_1}$
%normalized by $H(u_{1,k_2}\otimes u_{1,k_1})= \min\{k_1,k_2\}$ }
%\end{align*}
% is true. What do you think of it?)}

%(2)
%If $\A  = [\ell]$ and $\B = \emptyset$, then $SST_{\A/\B}(k_j) = B^{1, k_j}$ for $j=1,2$. In this case, we can check that the map ${\mathsf H}$ on $B^{1, k_1}\times B^{1, k_2}$ coincides with the local energy function $H$ on  %$B^{1, k_2}\otimes B^{1, k_1}$ (in reverse order) normalized by $H(u_{1,k_2}\otimes u_{1,k_1})=-\min\{k_1,k_2\}$ (cf.\ \cite[Section 3]{NaYa}). In both cases (1) and (2), $\upvarsigma$ is equal to the combinatorial $R$-matrix %$\sigma$.

(3) Suppose that $\A$ is finite with $|\A_0|=n$ and $|\A_1|=m$ and $\B=\emptyset$. Then $\SST_{\A/\B}(k)$ for $k\geq 1$ can be viewed as a crystal over the quantum superalgebra $U_q(\gl_{m|n})$ \cite{BKK}.  It would be very nice to find a representation theoretical meaning of ${\mathsf H}$ and $\upvarsigma$ from finite-dimensional modules over the quantum affine superalgebra $U_q(\widehat{\mf{sl}}_{m|n})$ \cite{Zh}.

(4) %Suppose that $\A=\Z_{>0}'$ and $\B=\Z_{<0}'$.  Recall that $\cP=SST_{\A/\B}(0)$ (see Section \ref{Sec:IrrChar} (1)).   {\color{red} (I think we can describe  ${\mathsf H}(\lambda,\mu)$ and $\upvarsigma(\lambda,\mu)$ for $\lambda,\mu\in \cP=SST_{\A/\B}(0)$ maybe in terms of Frobenius notations of partitions.... or give some concrete examples)}
Let ${\bf T}= (T_1, T_2, T_3)$ be as in Example \ref{Ex: RSK}. Since
\begin{align*}
\mathrm{wt}(T_2) &= \epsilon_{0'} + \epsilon_{ 2'} + \epsilon_{ 6'}  + \epsilon_{ 7'} - \epsilon_{\hyphen 4'} - \epsilon_{\hyphen 2'} - \epsilon_{\hyphen 1'}, \\
\mathrm{wt}(T_3) &= \epsilon_{ 4'}  + \epsilon_{ 5'} - \epsilon_{\hyphen 2'} - \epsilon_{\hyphen 1'},
\end{align*}
we have $\mathsf{s}_{T_2,T_3} = ( + \ + \ - \ - \ + \ + \ - \ + \ - \ + \ -  )$. Thus, the reduced sequence
$\mathsf{s}^{\rm red} = (\, + \, )$ gives $\mathsf{H}(T_2, T_3) = 0$ and
$$ \sigma(T_2, T_3 ) = (T_2', T_3') \in \SST_{\A/ \B}(0) \times \SST_{\A/\B}(1) ,$$
where
\begin{align*}
T_2' =
 (\
\resizebox{.10\hsize}{!}{${\def\lr#1{\multicolumn{1}{|@{\hspace{.6ex}}c@{\hspace{.6ex}}|}{\raisebox{-.1ex}{$#1$}}}\raisebox{ .1ex}
{$\begin{array}{ccc}
\cline{1-3}
\lr{0'}  & \lr{ 6' } & \lr{ 7' } \\
\cline{1-3}
\end{array}$}} $}
\ , \
\resizebox{.125\hsize}{!}{${\def\lr#1{\multicolumn{1}{|@{\hspace{.6ex}}c@{\hspace{.6ex}}|}{\raisebox{-.1ex}{$#1$}}}\raisebox{ .1ex}
{$\begin{array}{|ccc}
\cline{1-3}
\lr{ \hyphen 4' }  & \lr{ \hyphen 2' } & \lr{ \hyphen 1'}  \\
\cline{1-3}
\end{array}$}} $}
\  )
, \quad
T_3' = (\
\resizebox{.1\hsize}{!}{${\def\lr#1{\multicolumn{1}{|@{\hspace{.6ex}}c@{\hspace{.6ex}}|}{\raisebox{-.1ex}{$#1$}}}\raisebox{ .1ex}
{$\begin{array}{ccc}
\cline{1-3}
\lr{2'} &  \lr{ 4' } & \lr{ 5' } \\
\cline{1-3}
\end{array}$}} $}
\ , \
\resizebox{.085\hsize}{!}{${\def\lr#1{\multicolumn{1}{|@{\hspace{.6ex}}c@{\hspace{.6ex}}|}{\raisebox{-.1ex}{$#1$}}}\raisebox{ .1ex}
{$\begin{array}{|cc}
\cline{1-2}
  \lr{ \hyphen 2' } & \lr{ \hyphen 1' } \\
\cline{1-2}
\end{array}$}} $}
\  ).
\end{align*}
%Note that the corresponding Young diagrams in the level 1 fermionic Fock space are $\lambda_{T_2'} = (8,8,3,1)$ and $\lambda_{T_3'} = (5,5,4)$.

In the same manner, we compute
\begin{align*}
\mathsf{s}_{T_1,T_2} &= (-\ -\ +\ +\ +\ -\ +\ +\ -\ +\ -\ -\ +\ +),\\
\mathsf{s}_{T_1,T_2'} &= (-\ -\  +\ +\ +\  +\ +\ -\ + \ -\ -\ +\ +),
\end{align*}
which yield $\mathsf{H}(T_1, T_2) = -2$ and $\mathsf{H}(T_1, T_2') = -2$.

While $\SST_{\A/\B}(k)$ for $k\in \Z$ produces a character of a level one integrable highest weight module over $U_q(\gl_\infty)$, it also corresponds to a KR module over $U_q(\widehat{\mf{sl}}_\infty)$, where $\widehat{\mf{sl}}_\infty$ is an affinization of $\mf{sl}_\infty$ \cite{Her}. As in (3), we expect that $\mathsf{H}$ and $\upvarsigma$ are closely related with the theory of KR modules over $U_q(\widehat{\mf{sl}}_\infty)$.

}
\end{ex}

Now, we fix $n\geq 1$.  For simplicity, we put
\begin{equation}
\begin{split}
\F^n&=\F_{\A/\B}^n,\\
\F^\mu &= \SST_{\A/\B}(\mu_1)\times \cdots \times \SST_{\A/\B}(\mu_n),
\end{split}
\end{equation}
for $\mu=(\mu_1,\ldots,\mu_n)\in \Z^n$.
Clearly, we have $\F^n=\bigsqcup_{\mu\in \Z^n}\F^{\mu}$.
For $1 \le i \le n-1$, let $\upvarsigma_i$ be the map on $\F^\mu\subset \F^{n}$, which acts as $\upvarsigma$ on $\SST_{\A/\B}(\mu_i) \times \SST_{\A/\B}(\mu_{i+1})$ and as identity elsewhere,
and let ${\mathsf H}_i$ be the map on $\F^n$ given by ${\mathsf H}_i(T_1, \ldots, T_n) ={\mathsf H}(T_i, T_{i+1})$ for $(T_1,\ldots, T_n)\in \F^n$.

We define a {\it combinatorial energy function} ${\mathsf D} :  \F^n \longrightarrow \Z$ by
$$
{\mathsf D} (\mathbf{T}) = \sum_{1\leq i<j\leq n}{\mathsf H}_{i}(\upvarsigma_{i+1}\upvarsigma_{i+2}\cdots \upvarsigma_{j-1}(\mathbf{T})) \ \ \ \ (\mathbf{T} \in \F^n).
$$
Then we have the following, which is a generalization of \cite{NaYa}. The proof is given in the next section.

\begin{thm}\label{Thm:Energy=Charge} For ${\bf T}\in \F^n$, we have
\begin{equation*}
{\mathsf D}({\bf T}) = - \mathrm{charge} (Q_{\bf T}),
\end{equation*}
where $Q_{\bf T}$ is the rational semistandard tableau corresponding to ${\bf T}$ under the RSK map $\kappa_{\A/\B}$ on $\F^n$ in Theorem \ref{Thm: RSK}. In particular, we have ${\mathsf D}({\bf T})={\mathsf D}({\bf T}')$ for ${\bf T}, {\bf T}'\in \F^n$ such that $Q_{\bf T}=Q_{\bf T'}$.
\end{thm}

\begin{ex} \label{Ex: energy2}
{\rm
Continuing Example \ref{Ex: energy} (4), we have
$$ \mathsf{D}(\mathbf{T}) = \mathsf{H}(T_1, T_2) + \mathsf{H}(T_1, T_2') + \mathsf{H}(T_2, T_3) = -4. $$
Since
$$\mathrm{charge}(Q_\mathbf{T}) =
c\left(\
\resizebox{.25\hsize}{!}{${\def\lr#1{\multicolumn{1}{|@{\hspace{.6ex}}c@{\hspace{.6ex}}|}{\raisebox{-.3ex}{$#1$}}}\raisebox{-.1ex}
{$\begin{array}{cccccccc}
\cline{1-7}
  \lr{\,1\,}  & \lr{\,1\,}  & \lr{\,1\,}  & \lr{\,1\,}  & \lr{\,1\,}  & \lr{\,2\,} & \lr{\,2 \,} \\
\cline{1-7}
 \lr{\,2\,}  & \lr{\,3\,}  & \lr{\,3\,}  \\
\cline{1-3}
\end{array}$}}$} \
\right) = 4,
$$
we have $\mathsf{D}(\mathbf{T})= - \mathrm{charge}(Q_\mathbf{T}) $.
}
\end{ex}

As a consequence, we obtain a combinatorial realization of \eqref{Eq:Q function} in terms of ${\mathsf D}$.

\begin{thm}\label{Thm:Q=gradedchar} For $\mu\in \Z_+^n$, we have
\begin{equation*}
Q^{\A/\B}_{\mu}= \sum_{ \mathbf{T} \in \F^\mu}q^{-{\mathsf D}( \mathbf{T} )}{\bf x}_{\A/\B}^{ \mathbf{T} },
\end{equation*}
where ${\bf x}_{\A/\B}^{\mathbf{T}} = \prod_{i=1}^n {\bf x}_{\A/\B}^{T_i}$ for $ \mathbf{T}=(T_1,\ldots,T_n)\in \F^\mu$.
\end{thm}
\pf
Restricting $\kappa_{\A/\B}$ to $\F^\mu$, we have a weight preserving   bijection:
$$ \kappa_{\A/\B} :  \F^\mu \longrightarrow  \bigsqcup_{\lambda \in \cP_{\A/\B, n}} \SST_{\A/\B}(\lambda) \times \SST_{[n]}(\lambda)_\mu .$$
Thus, the assertion follows from $\eqref{Eq:Q function}$, $\eqref{Eq: KF poly}$, and Theorem \ref{Thm:Energy=Charge}.
\qed

\subsection{Proof of Theorem \ref{Thm:Energy=Charge}}

The proof is given in two steps. We will first define a regular $A_{n-1}$-crystal structure on $\F^n$ and show that $-{\mathsf D}$ is equal to the charge on $\F^n$ as an $A_{n-1}$-crystal.
Next, we will show that the RSK type correspondence $\kappa_{\A/\B}$ for parabolically semistandard tableaux is an $A_{n-1}$-crystal isomorphism, which is a key part in the proof of Theorem \ref{Thm:Energy=Charge}.
%This isomorphism allows us to consider Theorem $\ref{Thm:Energy=Charge}$ from a viewpoint of crystal bases theory.

Let us define an $A_{n-1}$-crystal structure on $\F^n$.
Let $\cM_{\A/\B, n}$ be the set of matrices $\mathbf{m}= (m_{ij})$ with non-negative integral entries $(i\in \A \sqcup \B,\ j \in [n])$ satisfying
(1) $\sum_{i,j} m_{ij} < \infty$, (2) $m_{ij} \in \{0,1\}$ if $i$ is odd,. Note that for ${\bf T}=(T_1,\ldots,T_n)\in \F^n$ with
\begin{equation*}
{\rm wt}_{\A/\B}(T_j)=\sum_{a\in \A}m_{aj}\epsilon_a - \sum_{b\in \B}m_{bj}\epsilon_b,
\end{equation*}
for $1\leq j\leq n$,  the map sending ${\bf T}$ to ${\bf m}=(m_{ij})$ gives a natural bijection from $\F^n$ to $\cM_{\A/\B, n}$.

Let ${\bf m}\in \cM_{\A/\B, n}$ be given.
For $i \in \A \sqcup \B$, let $ \mathbf{m}_{(i)} = ( m_{i j}  )_{j\in [n]}$ be the $i$th row of ${\bf m}$, and set $|\mathbf{m}_{(i)}| = \sum_{j\in [n]} m_{i j}$. Let $\lambda^{(i)}\in \Z_+^n$ be given by
$$
\lambda^{(i)} = \left\{
              \begin{array}{ll}
                (|\mathbf{m}_{(i)}|,0,\ldots,0),  & \hbox{ if $i\in \A_0$}, \\
                (1^{|\mathbf{m}_{(i)}|},0,\ldots,0),  & \hbox{ if $i\in \A_1$},\\
                (0,\ldots,0,-|\mathbf{m}_{(i)}|),  & \hbox{ if $i\in \B_0$}, \\
                (0,\ldots,0,-1^{|\mathbf{m}_{(i)}|}),  & \hbox{ if $i\in \B_1$}.
              \end{array}
            \right.
$$
We identify ${\bf m}_{(i)}$ with a unique rational semistandard tableau $T^{(i)}\in \SST_{[n]}(\lambda^{(i)})$ such that
$$
{\rm wt}_{[n]}\left(T^{(i)}\right) = \left\{
              \begin{array}{ll}
               \ \ \,\sum_{j\in [n]}m_{ij}\epsilon_j,  & \hbox{ if $i\in \A$}, \\
               -\sum_{j\in [n]}m_{ij}\epsilon_j,  & \hbox{ if $i\in \B$}.
              \end{array}
            \right.
$$

For example, if $a\in \A_0$ and $b\in \B_1$,  then we have

\begin{equation*}
\begin{split}
& \mathbf{m}_{(a)} = (2,0,1,2) \ \ \longleftrightarrow \ \ T^{(a)} = \ \
\resizebox{.17\hsize}{!}{${\def\lr#1{\multicolumn{1}{|@{\hspace{.6ex}}c@{\hspace{.6ex}}|}{\raisebox{-.1ex}{$#1$}}}\raisebox{ .1ex}
{$\begin{array}{|ccccc}
\cline{1-5}
\lr{\, 1\,}  & \lr{\, 1\,} & \lr{\, 3\,} & \lr{\, 4\,} & \lr{\, 4\,}  \\
\cline{1-5}
\end{array}$}} $}\ \
\in \SST_{[4]} (5,0,0,0), \\
& \mathbf{m}_{(b)}= (1,1,0,1)  \ \ \longleftrightarrow \ \ T^{(b)} = \ \
\resizebox{.055\hsize}{!}{${\def\lr#1{\multicolumn{1}{|@{\hspace{.6ex}}c@{\hspace{.6ex}}|}{\raisebox{-.3ex}{$#1$}}}\raisebox{-.6ex}
{$\begin{array}{c|c}
\cline{1-1}
\lr{\hyphen 4} &     \\
\cline{1-1}
\lr{\hyphen 2} &     \\
\cline{1-1}
\lr{\hyphen 1}&     \\
\cline{1-1}
\end{array}$}}$}
 \in \SST_{[4]} (0,\hyphen1, \hyphen1,\hyphen1).
\end{split}
\end{equation*}
Then we define a regular $A_{n-1}$-crystal structure on $\cM_{\A/\B, n}$ and hence on $\F^n$ via the correspondence
\begin{align} \label{Eq: matrix crystal}
\mathbf{m}\ \  \longleftrightarrow  \ \
\overset{\longleftarrow}{\bigotimes}_{a\in \A} T^{(a)}\ \otimes \
\overset{\longrightarrow}{\bigotimes}_{b\in \B} T^{(b)}. %\  \in  \ \
%\overset{\longleftarrow}{\bigotimes_{a\in \A }} SST_{[n]}(\lambda^{(a)}) \otimes \overset{\longrightarrow}{\bigotimes_{b\in \B }} SST_{[n]}(\lambda^{(b)}),
\end{align}
Here we understand $\overset{\longleftarrow}{\bigotimes}_{a \in \A}T^{(a)}$ as a tenor product with respect to the reverse total order on $\A$. Since $T^{(a)}$ is an empty tableau except for finitely many
$a\in \A$, it is well-defined. Similarly, $\overset{\longrightarrow}{\bigotimes}_{b \in \B}T^{(b)}$ is a tensor product with respect to the total order on $\B$. One may assume that the row indices of ${\bf m}\in {\cM_{\A/\B,n}}$ are parametrized by $\B^\pi \ast \A$, and we read each row in ${\bf m}$ from bottom to top.

Since $ \mathcal{F}^n$ is an $A_{n-1}$-crystal,
one can consider $ \upvarepsilon_\alpha (\mathbf{T})$ and $\upvarphi_\alpha (\mathbf{T})$ for $\mathbf{T} \in \mathcal{F}^n$ and $\alpha \in \Delta^+_{n-1}$  as in $\eqref{Eq: epsilon and phi}$.
By definitions of ${\mathsf H}$ and $\upvarsigma$, we can check that
\begin{align} \label{Eq: equalities}
{\mathsf H}_{i}(\mathbf{T}) = -\min\{\upvarepsilon_i(\mathbf{T}),\upvarphi_i(\mathbf{T})\}, \qquad \upvarsigma_i(\mathbf{T}) = S_i (\mathbf{T}),
\end{align}
for $1\leq i\leq  n-1$. In particular $\upvarsigma_i$'s satisfy the braid relations.
Thus, combining $\eqref{Eq: equalities}$ with $ \upvarepsilon_\alpha (\mathbf{T})$ and $\upvarphi_\alpha (\mathbf{T})$,
we obtain the following, which generalizes Proposition \ref{prop:energy description}.

\begin{prop}
For $ \mathbf{T} \in \mathcal{F}^n$, we have
\begin{equation*}
{\mathsf D} (\mathbf{T}) =-\sum_{\alpha\in \Delta^+_{n-1}}\min\{ \upvarepsilon_\alpha(\mathbf{T}), \upvarphi_\alpha(\mathbf{T})\}.
\end{equation*}
\end{prop}

Moreover, the regular $A_{n-1}$-crystal structure of $ \mathcal{F}^n$ enables us to consider the charge of $ \mathbf{T}\in \F^n$.
By Theorem \ref{Thm: charge}, we have
\begin{cor} \label{Cor:charge}
For $ \mathbf{T} \in \mathcal{F}^n$, we have
${\mathsf D} (\mathbf{T})=- {\rm charge}( \mathbf{T} ) $.
\end{cor}
This also immediately implies that ${\mathsf D}\circ \upvarsigma_i = {\mathsf D} $ for $1\leq i\leq n-1$.

Next, we interpret the map $\kappa_{\A/\B}$
$$ \kappa_{\A/\B} :  \mathcal{F}^n \longrightarrow  \bigsqcup_{\lambda \in \cP_{\A/\B, n}} \SST_{\A/\B}(\lambda) \times \SST_{[n]}(\lambda). $$
from a viewpoint of crystal.
We assume that the right-hand side is an $A_{n-1}$-crystal, where the operators $\tte_j$ and $\ttf_j$ act on the second factor $\SST_{[n]}(\lambda)$.

\begin{thm} \label{Thm:crystal isom kappa}
The map $\kappa_{\A/\B}$ is an $A_{n-1}$-crystal isomorphism.
\end{thm}
\pf  Let us recall the bijections $\rhocol$ and $\rhorow$ given in \eqref{Eq:rho}. Suppose $\B = \emptyset$ and write $\mathcal{F}_{\A} = \mathcal{F}_{\A/\emptyset}$.
Extending $\rhocol$ and $\rhorow$ to $\mathcal{F}_{\A}^n$, we have bijections
$$
\rhocol, \rhorow : \mathcal{F}_{\A}^n \longrightarrow  \bigsqcup_{\lambda \in \cP_{\A}\cap\cP_n} \SST_{\A}(\lambda) \times \SST_{[n]}(\lambda),
$$
which are indeed the RSK correspondences since $\F^n_{\A}$ can be identified with $\cM_{\A/\emptyset, n}$. Let ${\bf T}\in \F^n_\A$ be given with the corresponding matrix ${\bf m}\in \cM_{\A/\emptyset,n}$.  If we write $\rhocol(\mathbf{T}) = (P_c, Q_c)$ and $\rhorow(\mathbf{T}) = (P_r, Q_r)$, then
\begin{align} \label{Eq: RSK1}
\overset{\longrightarrow}{\bigotimes}_{a\in \A} {\bf m}_a \equiv Q_c, \qquad \overset{\longleftarrow}{\bigotimes}_{a\in \A} {\bf m}_a \equiv Q_r.
\end{align}
Indeed, the first equivalence follows from \cite[Theorem 3.11]{K07} on $(\gl_{m|n},\gl_{u|v})$-bicrystal isomorphism over general linear Lie superalgebras, where we replace $\gl_{u|v}$ with $\gl_{n|0}$ and $\gl_{m|n}$ with a Lie superalgebra $\gl_{\A}$ associated to $\A$ (see also \cite[Lemma 4.9]{K09}). Similarly, the second equivalence can be obtained by changing the parity of $\A$ and applying \cite[Theorem 4.5]{K07}.

Let $\mathbf{T} \in \mathcal{F}^n$ be given, and ${\bf m}$ the corresponding matrix in $\cM_{\A/\B,n}$.
Let $\kappa_{\A/\B}(\mathbf{T}) = (P_{\mathbf{T}}, Q_{\mathbf{T}})$.
We keep the same notations $Q$, $Q^{\vee}$, and $U_R$ in ($\kappa$-1)-($\kappa$-4) in Section \ref{Sec:PSST}.
It follows from $\eqref{Eq:sigma delta}$, ($\kappa$-1), ($\kappa$-2), and the first equivalence in $\eqref{Eq: RSK1}$ that
$$ Q^{\vee}
\equiv \left( \overset{\longrightarrow}{\bigotimes}_{b\in \B^{\pi}}  \mathbf{m}_b^\vee  \right)^{\vee}
\equiv \left( \overset{\longleftarrow}{\bigotimes}_{b\in \B}  \mathbf{m}_b^\vee \right)^{\vee}
\equiv  \overset{\longrightarrow}{\bigotimes}_{b\in \B} \mathbf{m}_b.$$
We should remark that on $\cM_{\emptyset/\B,n}$ the $A_{n-1}$-crystal structure is dual to that of $\cM_{\A/\emptyset,n}$.
Moreover, since $k < a$ in $[n]*\A$ for all $k\in [n]$ and $a\in \A$, $\eqref{Eq:sigma delta}$ and the second equivalence in $\eqref{Eq: RSK1}$ give
$$
Q_{\mathbf{T}}
\equiv U_R
\equiv \left(\overset{\longleftarrow}{\bigotimes}_{a\in \A}  \mathbf{m}_a \right)\otimes Q^\vee
\equiv \left(\overset{\longleftarrow}{\bigotimes}_{a\in \A}  \mathbf{m}_a\right) \otimes
\left(\overset{\longrightarrow}{\bigotimes}_{b\in \B} \mathbf{m}_b \right)
\equiv \mathbf{T},
$$
which implies that $\kappa_{\A/\B}$ is a morphism of $A_{n-1}$-crystals. Thus the assertion follows from Theorem \ref{Thm: RSK}.
\qed

\begin{cor}\label{Cor:charge-2}
For ${\bf T}\in \F^n$,  $\mathrm{charge}({\bf T})=\mathrm{charge}(Q_{\bf T})$.
\end{cor}

Now, Theorem \ref{Thm:Energy=Charge} follows from Corollaries \ref{Cor:charge} and \ref{Cor:charge-2}. This completes the proof.

\begin{ex}
{\rm
Continuing Example \ref{Ex: RSK}, we have the matrix $\mathbf{m}$ corresponding to $\mathbf{T}$ as follows:
\begin{align*}
\mathbf{m} = \begin{array}{c|ccc}
                & \ 1\  &\  2\  & \ 3\ \\ \hline
               \hyphen 4'\ &\ \bullet\ &\ \bullet\ &\ \cdot \\
               \hyphen 3'& \ \bullet\ &\ \cdot\ &\ \cdot \\
               \hyphen 2' &\ \cdot\ &\ \bullet\ &\ \bullet \\
               \hyphen 1' &\ \cdot\ &\ \bullet\ &\ \bullet \\
               0' & \ \bullet\ &\ \bullet\ &\ \cdot \\
               1' & \ \bullet\ &\ \cdot\ &\ \cdot \\
               2' & \ \cdot\ &\ \bullet\ &\ \cdot \\
               3' & \ \bullet\ &\ \cdot\ &\ \cdot \\
               4' & \ \bullet\ &\ \cdot\ &\ \bullet \\
               5' & \ \bullet\ &\ \cdot\ &\ \bullet \\
               6' & \ \cdot\ &\ \bullet\ &\ \cdot \\
               7' & \ \cdot\ &\ \bullet\ &\ \cdot
             \end{array}
\end{align*}
\vskip 2mm
\noindent where $\cdot$ and $\bullet$ denote $0$ and $1$ respectively. This yields
\begin{align*}
&  T^{(0')} =
\resizebox{.05\hsize}{!}{${\def\lr#1{\multicolumn{1}{|@{\hspace{.6ex}}c@{\hspace{.6ex}}|}{\raisebox{-.3ex}{$#1$}}}\raisebox{-.0ex}
{$\begin{array}{c}
\cline{1-1}
\lr{\,1\,}   \\
\cline{1-1}
\lr{\,2\,}   \\
\cline{1-1}
\end{array}\ ,$}}$}\ \
T^{(1')} = T^{(3')} =
\resizebox{.05\hsize}{!}{${\def\lr#1{\multicolumn{1}{|@{\hspace{.6ex}}c@{\hspace{.6ex}}|}{\raisebox{-.3ex}{$#1$}}}\raisebox{-.0ex}
{$\begin{array}{c}
\cline{1-1}
\lr{\,1\,}   \\
\cline{1-1}
\end{array}\ ,$}}$}\ \
T^{(2')} = T^{(6')} = T^{(7')} =
\resizebox{.05\hsize}{!}{${\def\lr#1{\multicolumn{1}{|@{\hspace{.6ex}}c@{\hspace{.6ex}}|}{\raisebox{-.3ex}{$#1$}}}\raisebox{-.0ex}
{$\begin{array}{c}
\cline{1-1}
\lr{\,2\,}   \\
\cline{1-1}
\end{array}\ ,$}}$}\ \
T^{(4')} = T^{(5')} =
\resizebox{.05\hsize}{!}{${\def\lr#1{\multicolumn{1}{|@{\hspace{.6ex}}c@{\hspace{.6ex}}|}{\raisebox{-.3ex}{$#1$}}}\raisebox{-.0ex}
{$\begin{array}{c}
\cline{1-1}
\lr{\,1\,}   \\
\cline{1-1}
\lr{\,3\,}   \\
\cline{1-1}
\end{array}\ ,$}}$}   \\
&  T^{(\hyphen 1')} =  T^{(\hyphen 2')} =
\resizebox{.05\hsize}{!}{${\def\lr#1{\multicolumn{1}{|@{\hspace{.6ex}}c@{\hspace{.6ex}}|}{\raisebox{-.3ex}{$#1$}}}\raisebox{-.0ex}
{$\begin{array}{c}
\cline{1-1}
\lr{\hyphen 3}   \\
\cline{1-1}
\lr{\hyphen 2}   \\
\cline{1-1}
\end{array}\ ,$}}$}\ \
T^{(\hyphen 3')} =
\resizebox{.05\hsize}{!}{${\def\lr#1{\multicolumn{1}{|@{\hspace{.6ex}}c@{\hspace{.6ex}}|}{\raisebox{-.3ex}{$#1$}}}\raisebox{-.0ex}
{$\begin{array}{c}
\cline{1-1}
\lr{\hyphen 1}   \\
\cline{1-1}
\end{array}\ ,$}}$}\ \
 T^{(\hyphen 4')} =
\resizebox{.05\hsize}{!}{${\def\lr#1{\multicolumn{1}{|@{\hspace{.6ex}}c@{\hspace{.6ex}}|}{\raisebox{-.3ex}{$#1$}}}\raisebox{-.0ex}
{$\begin{array}{c}
\cline{1-1}
\lr{\hyphen 2}   \\
\cline{1-1}
\lr{\hyphen 1}   \\
\cline{1-1}
\end{array}\ .$}}$}
\end{align*}
Thus, using the Schensted's bumping algorithm, we see
\begin{align*}
& T^{(7')} \otimes T^{(6')} \otimes T^{(5')} \otimes T^{(4')} \otimes T^{(3')} \otimes T^{(2')} \otimes T^{(1')} \otimes T^{(0')}
 \equiv
\resizebox{.25\hsize}{!}{${\def\lr#1{\multicolumn{1}{|@{\hspace{.6ex}}c@{\hspace{.6ex}}|}{\raisebox{-.3ex}{$#1$}}}\raisebox{-.1ex}
{$\begin{array}{cccccccc}
\cline{1-7}
\lr{\,1\,}  & \lr{\,1\,}  & \lr{\,1\,}  & \lr{\,1\,}  &  \lr{\,1\,} & \lr{\,2\,} & \lr{\,2\,}  \\
\cline{1-7}
 \lr{\,2\,}  &  \lr{\,2\,} & \lr{\,3\,} & \lr{\,3\,}    \\
\cline{1-4}
\end{array}$}}$} \  ,
\\
& T^{ ( \hyphen 4')} \otimes T^{(\hyphen 3')} \otimes T^{(\hyphen 2')} \otimes T^{(\hyphen 1')} \equiv
\resizebox{.03\hsize}{!}{${\def\lr#1{\multicolumn{1}{|@{\hspace{.6ex}}c@{\hspace{.6ex}}|}{\raisebox{-.3ex}{$#1$}}}\raisebox{-.0ex}
{$\begin{array}{c}
\cline{1-1}
\lr{ \,3\,}   \\
\cline{1-1}
\end{array}$}}$}
\otimes
\resizebox{.03\hsize}{!}{${\def\lr#1{\multicolumn{1}{|@{\hspace{.6ex}}c@{\hspace{.6ex}}|}{\raisebox{-.3ex}{$#1$}}}\raisebox{-.0ex}
{$\begin{array}{c}
\cline{1-1}
\lr{ \,2\,}   \\
\cline{1-1}
\lr{ \,3\,}   \\
\cline{1-1}
\end{array}$}}$}
\otimes
\resizebox{.03\hsize}{!}{${\def\lr#1{\multicolumn{1}{|@{\hspace{.6ex}}c@{\hspace{.6ex}}|}{\raisebox{-.3ex}{$#1$}}}\raisebox{-.0ex}
{$\begin{array}{c}
\cline{1-1}
\lr{ \,1\,}   \\
\cline{1-1}
\end{array}$}}$}
\otimes
\resizebox{.03\hsize}{!}{${\def\lr#1{\multicolumn{1}{|@{\hspace{.6ex}}c@{\hspace{.6ex}}|}{\raisebox{-.3ex}{$#1$}}}\raisebox{-.0ex}
{$\begin{array}{c}
\cline{1-1}
\lr{\,1\,}   \\
\cline{1-1}
\end{array}$}}$}
\equiv
\resizebox{.13\hsize}{!}{${\def\lr#1{\multicolumn{1}{|@{\hspace{.6ex}}c@{\hspace{.6ex}}|}{\raisebox{-.3ex}{$#1$}}}\raisebox{-.1ex}
{$\begin{array}{ccccc}
\cline{1-4}
\lr{\,1\,}  & \lr{\,1\,}  & \lr{\,2\,}  & \lr{\,3\,}   \\
\cline{1-4}
 \lr{\,3\,}   \\
\cline{1-1}
\end{array}$}}$}
\equiv  Q^{\vee} \ .
\end{align*}
Therefore, we have
\begin{align*}
\mathbf{T} &\equiv \overset{\longleftarrow}{\bigotimes}_{a\in \Z_{\ge0}'} T^{(a)}\ \otimes \
\overset{\longrightarrow}{\bigotimes}_{b\in \Z_{<0}'} T^{(b)} \\
&\equiv
\resizebox{.25\hsize}{!}{${\def\lr#1{\multicolumn{1}{|@{\hspace{.6ex}}c@{\hspace{.6ex}}|}{\raisebox{-.3ex}{$#1$}}}\raisebox{-.1ex}
{$\begin{array}{cccccccc}
\cline{1-7}
\lr{\,1\,}  & \lr{\,1\,}  & \lr{\,1\,}  & \lr{\,1\,}  &  \lr{\,1\,} & \lr{\,2\,} & \lr{\,2\,}  \\
\cline{1-7}
 \lr{\,2\,}  &  \lr{\,2\,} & \lr{\,3\,} & \lr{\,3\,}    \\
\cline{1-4}
\end{array}$ }} $}   \otimes \
\resizebox{.13\hsize}{!}{${\def\lr#1{\multicolumn{1}{|@{\hspace{.6ex}}c@{\hspace{.6ex}}|}{\raisebox{-.3ex}{$#1$}}}\raisebox{-.1ex}
{$\begin{array}{ccccc}
\cline{1-4}
\lr{\,1\,}  & \lr{\,1\,}  & \lr{\,2\,}  & \lr{\,3\,}   \\
\cline{1-4}
 \lr{\,3\,}   \\
\cline{1-1}
\end{array}$}}$}
\\
&\equiv Q_{\mathbf{T}}\ .
\end{align*}
}
\end{ex}

\end{document}